\documentclass[11pt, final, reqno, letterpaper]{article}
\usepackage[english]{babel}
\usepackage{amsmath,amssymb,dsfont, amsthm}
\usepackage{epsfig}
\usepackage{bbm}
\usepackage{natbib}
\usepackage{color}
\usepackage{multirow}
\usepackage{color}

\graphicspath{{Figures/}}

\definecolor{purple}{rgb}{0.55,0.2,0.90}

\let\proglang=\textsf
\newcommand{\pkg}[1]{{\fontseries{b}\selectfont #1}}

\renewcommand{\hat}{\widehat}

\newcommand{\LL}{{\mathbb L}}

\newcommand{\RR}{\mathbb{R}}
\newcommand{\UU}{\mathbb{U}}

\newcommand{\YY}{{\mathbb Y}}
\newcommand{\WW}{{\mathbb W}}
\newcommand{\HH}{{\mathbb H}}
\newcommand{\NN}{\mathbb{N}}

\newtheorem{theorem}{Theorem}[section]
\newtheorem{lemma}[theorem]{Lemma}
\newtheorem{prop}[theorem]{Proposition}

\title{On asymptotics of the discrete convex LSE of a pmf}
\author{Fadoua Balabdaoui$^1$, C\'ecile Durot$^2$\thanks{Corresponding author. Email address: cecile.durot@gmail.com} and Fran\c{cois} Koladjo$^{3,4}$}

\begin{document}
\maketitle

\begin{center}
$^{1}$CEREMADE, Universit\'e Paris-Dauphine, 75775, Paris, France \\
$^{2}$UFR SEGMI, Universit\'e Paris Ouest Nanterre La D\'efense, F-92001, Nanterre, France\\
$^{3}$Equipe Probabilit\'e, Statistique et Mod\'elisation, UMR CNRS 8628, Universit\'e de Paris-Sud, 91405 Orsay Cedex, France \\
$^{4}$CIPMA-Chaire UNESCO, FAST, UAC, 072BP50 Cotonou, B\'enin.\\
\end{center}

\begin{abstract}
In this article, we derive the weak limiting distribution of the least squares estimator (LSE) of a convex probability mass function (pmf) with a finite support. We show that it can be defined via a certain convex projection of a Gaussian vector. Furthermore, samples of any given size from this limit distribution can be generated using an efficient Dykstra-like algorithm. 
\end{abstract}

\bigskip

Keywords:  convex, least squares, nonparametric estimation, pmf, shape-constraints

\section{Introduction}
Non-parametric estimation under a shape constraint of a density on a given sub-interval of $\RR$, has attracted considerable attention over the past decades. Typical shape constraints are monotonicity, convexity, log-concavity. Typical estimators are the maximum likelihood estimator (MLE) and the least-squares estimator (LSE). Both of them are obtained by minimization of a given criterion over the set of all densities that satisfy the considered shape constraint. Even if the MLE and LSE  uniquely exist, no closed form is available for these estimators so a key step is to provide a precise characterization of the estimators as well as an algorithm for practical implementation. \cite{grenander1956theory} first gives such a characterization for the MLE of a monotone density, and the pointwise weak convergence of the MLE is derived in \cite{rao1969estimation}. The estimator can easily been implemented using the Pool Adjacent Violators Algorithm as described in \cite{barlow1972statistical}. Both characterization and pointwise weak convergence of a convex density on the half-real line are investigated in \cite{gjw01B}, and practical implementation is discussed in \cite{groeneboom2008support}. The MLE of a log-concave density is characterized in \cite{dumbgen2009maximum} while its pointwise weak convergence is studied in \cite{balabdaoui2009limit}. Algorithmic aspects are treated in \cite{dumbgen2010logcondens}. In the aforementioned continuous case of estimating a density under a shape constraint over a given sub-interval of $\RR$, the global asymptotic behavior of estimators has also attracted attention. Precisely, the limit behavior of a distance between the estimator and the true density has been investigated. We refer to \cite{kulikov2005asymptotic} and  \cite{durot2012limit} for the limit distribution of the $L_{p}$-distance and the supremum-distance respectively, in the case of a monotone density. The rate of convergence of the supremum-distance between the log-concave MLE and the true density is given in \cite{dumbgen2009maximum}.

\medskip
More recently, attention has been given to estimation of a discrete probability mass function (pmf) under a shape constraint. Similar to the continuous case, no closed form is available for shape constrained estimators of a pmf so one needs a precise characterization. In contrast to the continuous case, the natural way to investigate the global limit behavior of the estimator is to compute the limit distribution of the whole process $\hat p_{n}-p_{0}$, where $\hat p_{n}$ is the considered estimator and $p_{0}$ is the true pmf. This approach was first considered by \cite{jankowski_09} in the case of a monotone pmf on $\NN$, and by \cite{BJRP2013} in the case of a log-concave pmf.
It should be noticed that the discrete case totally differs from the continuous case. In particular, the rate of convergence  is typically $\sqrt n$ (where $n$ denotes the sample size) in the discrete case whereas it is of smaller order in the continuous case. The characterization and the rate of convergence of the LSE of a convex pmf are given in \cite{durotetal_13} together with an algorithm for implementing the LSE, but the limit distribution of the LSE remains unknown.

\medskip
The aim of this paper is to compute the limit distribution, as the sample size grows to infinity, of the least-squares estimator (LSE) of a convex probability mass function (pmf). To be more precise, let $X_1, \ldots, X_n$ be i.i.d. from an unknown  discrete pmf $p_0$ whose support  takes the form $ \NN \cap [\kappa, \infty)  = \{\kappa, \kappa +1, \ldots, \}$ or $\NN \cap [\kappa, S]  = \{\kappa, \kappa+1, \ldots, S\}$   for some integers $S > \kappa \ge 0$. Here, $\kappa$ is assumed to be known whereas $S$ is unknown. 
Assuming  that $p_{0}$ is convex on  $ \NN \cap [\kappa, \infty)$, we are interested in the limiting behavior of 
the LSE of $p_0$.
The case $\kappa =1$ is of a particular interest in \cite{durotetal_14}, where  the problem of estimating the total number $N$ of species in a given area is investigated under a convexity constraint. Precisely, the distribution of  the abundance (that is, the number of  individuals) of a given species in a sample is assumed to be convex, so that the distribution $p^{+}$ of the zero-truncated counts (that is, the abundances of species that are present in the sample) is convex on  $ \NN \cap [1, \infty)$. Only zero-truncated counts are observable, so the authors consider the LSE of $p^+$ and, using a natural constraint that makes the problem identifiable, they derive two estimators as well as asymptotic confidence intervals for $N$. 

\medskip

Motivated by the aforementioned  application, we derive the weak limiting distribution of the LSE of $p_{0}$ under the convexity constraint. The limiting distribution is described as a piecewise convex approximation of a Gaussian process, where the pieces  are connected to the points of strict convexity of $p_{0}$. The Gaussian process involved in the limiting distribution depends on $p_{0}$ as well. For a given convex pmf $p_{0}$, we provide an algorithm for simulating the limiting distribution of the LSE. This amounts to simulate (many times)  the Gaussian process and its piecewise convex approximation. Since computing the piecewise convex approximation amounts to minimize the least-squares criterion over the intersection of closed convex cones, our algorithm combines two previous algorithms. The first one is implemented  in the package {\tt cobs} for R; see \cite{ng2007fast} for a full description. The function {\tt conreg} in the package {\tt cobs}  is used to minimize the least squares criterion over the closed convex cone of discrete functions that are convex on a given interval. Once an algorithm is available for minimizing the least squares criterion over a closed convex cone, the iterative algorithm by \cite{dykstra1983algorithm} is used to minimize the criterion over  the intersection of closed convex cones, thus providing the piecewise convex approximation. We use our algorithm to illustrate on a simulation study, the convergence of the distribution of the LSE to the limit distribution.

\medskip

The paper is organized as follows. All main notation are collected in Section \ref{sec:notation}. In Section \ref{sec:characterization}, we recall the characterization of the LSE obtained in \cite{durotetal_13} and we derive the $\sqrt n$-rate of convergence. In addition, we prove that any knot (that is, a point of strict convexity) of the true pmf is also (almost surely, for large enough $n$) a knot of the LSE. This allows us to characterize the support of the LSE in the case when the true pmf has a finite support. Section \ref{sec:asymptotics} is devoted to the weak convergence of the LSE. The limit distribution is computed in the general case and we investigate how the limit distribution simplifies in some specific cases, such as pmf having consecutive knots. Simulations are reported in Section \ref{sec:implementation}. We first investigate on a few examples whether the knots of the estimator include all true knots when the sample size is finite. Then, we illustrate the convergence of the distribution of the estimator to the limit distribution. All proofs are postponed to Section \ref{sec:proofs}.

\section{Notation}\label{sec:notation}

We feel it would benefit the reader if all main notation used in the paper were collected in a single section. 

\begin{description}
\item $\bullet$ \ In what follows, we will assume that observations $X_{1},\dots,X_{n}$ are i.i.d. from an unknown probability mass function (pmf) $p_{0}$ on $\NN$. We denote by $p_n$ the empirical pmf based on $X_{1},\dots,X_{n}$, that is
\begin{eqnarray}\label{eq:empirical}
p_n(j) = \frac{1}{n} \sum_{i=1}^n \mathbbm{1}_{\{X_i = j\}} ,\ j \in \NN
\end{eqnarray}
and by $\hat p_{n}$ the convex least squares estimator of $p_{0}$ (see Section \ref{sec:characterization} below for a precise definition).
\item $\bullet$ \ For an arbitrary pmf $p = \{p(k), k \in \NN\}$, we denote by $F_{p}$ the corresponding cumulative distribution function (cdf) given by
\begin{equation}\label{def:F_p}
F_{p}(k)=\sum_{j=0}^kp(j)
\end{equation}
for all $k\in\NN$ with $F_{p}(-1)=0$,
and we define
\begin{equation}\label{def:H_p}
H_{p}(z)=\sum_{k=0}^{z-1}F_{p}(k)
\end{equation}
for all $z\in\NN$ with the convention $H_{p}(0)=0$.
\item $\bullet$ \ For an arbitrary sequence $p = \{p(k), k \in \NN\}$, we denote by $\|p\|_{r}$ the $\ell_r$-norm of $p$, that is
\begin{eqnarray*}
\Vert p \Vert_r =
\begin{cases}
 \left(\sum_{k=0}^\infty \vert p(k) \vert^r \right)^{1/r} &  \text{if $r \in \NN \backslash \{0\}$} \\
  \sup_{k \in \NN} \vert p(k) \vert &  \text{if $r= \infty $}.
\end{cases}
\end{eqnarray*}

\item $\bullet$ \ For an arbitrary sequence $p = \{p(k), k \in \NN\}$ and $k \in  \NN\backslash\{0\}$, we denote by 
$$\Delta p(k)=p(k+1)-2p(k)+p(k-1)$$
the discrete Laplacian.  

\item $\bullet$ \ We denote by $\mathcal C$ the class of convex  sequences $p$ such that $\Vert p \Vert_2 < \infty$,  that is,
the set of sequences  $p$ satisfying 
\begin{equation}\label{condconvex}
\Vert p \Vert_2 < \infty \mbox{ and }\Delta p(k)\geq  0\mbox{ for all integers } k\geq 1. 
\end{equation}
It should be noted that any $p\in\cal C$ has to be non-negative and non-increasing. 


\item $ \bullet $ \ A point $k\geq 1$ in the support of $p \in \cal C$ is called a knot of $p$ if $\Delta p(k)>0$. 

\item $ \bullet $ \ $\UU = \UU(\omega)$ denotes a realization of a standard Brownian bridge from $(0,0)$ to $(1, 0)$ for $\omega \in \Omega$ a probability space. Define also the related Gaussian process $$\WW(k) = \UU(F_{p_0}(k)) - \UU(F_{p_0}(k-1)), k \in \NN,$$ with the convention $F_{p_{0}}(-1)=0$. 

\item $ \bullet $ \ Also, define
\begin{eqnarray}\label{eq:defH}
\HH(z) = \sum_{k=0}^{z-1} \UU(F_{p_0}(k)) = \sum_{k=0}^{z-1}\sum_{j=0}^{k} \WW(j)
\end{eqnarray}
with $\HH(0) =0$.  
\item $\bullet$ \ For all real numbers $x$, $x_{+}=\max\{x,0\}$.
\end{description}

\section{Characterization and tightness of the convex LSE}\label{sec:characterization}

Let $X_1, \ldots, X_n$ be i.i.d. from an unknown convex discrete pmf $p_0$ with support included in $\NN$. More precisely, we suppose that there exist some integers $S > \kappa \ge 0$  such that the support of $p_0$ takes the form $ \NN \cap [\kappa, \infty)  = \{\kappa, \kappa +1, \ldots, \}$ or $\NN \cap [\kappa, S]  = \{\kappa, \kappa+1, \ldots, S\}$, and that $p_{0}$ is convex on  $ \NN \cap [\kappa, \infty)$.  Here, $\kappa$ is assumed to be known whereas $S$ is unknown. The assumption $S > \kappa$ is made in order to avoid  the uninteresting situation of having to deal with a Dirac distribution. Since convexity is preserved under translation, we can assume without loss of generality that $\kappa =0$: in case $\kappa > 0$, the characterization as well as the asymptotic results for the LSE of the true convex pmf can be easily deduced from the ones established below using the simple fact that the support pmf of $X_i - \kappa$ admits 0 as its left endpoint. Thus, in the sequel we restrict our attention to the case of a convex pmf $p_{0}$ on $\NN$ with an unknown support.

Based on the sample $X_1, \ldots, X_n$, we consider the empirical pmf $p_n$ given in \eqref{eq:empirical}. We are mainly interested in the asymptotics of $\widehat p_n$, the least-squares estimator (LSE) of $p_0$ defined as the unique minimizer of the criterion 
\begin{eqnarray*}
\Phi_n(p) = \frac{1}{2} \sum_{j \in \NN} (p_n(j) - p(j))^2
\end{eqnarray*}
over $\mathcal C$, where we recall that $\mathcal C$ is the set of convex sequences $p$ such that $\Vert p \Vert_2 < \infty$. Existence and uniqueness of $\widehat p_n$ follows from the Hilbert projection Theorem,  see e.g. \cite{durotetal_13}, Section 2.1.

It has been proved in \cite{durotetal_13}, Theorem 1, that $\widehat p_n$ is a proper pmf in the sense that 
$$\sum_{j\geq 0}\hat p_{n}(j)=1.$$
This fact is very convenient because it means that in order to compute the estimator, we can minimize the criterion $\Phi_{n}$ over $\cal C$ rather than over the more constrained set of pmf's in $\cal C$. This allows us to use simpler algorithms.
It also has the advantage of giving more flexibility when deriving the characterizing Fenchel conditions for $\widehat p_n$. As it is the case in many shape constrained problems, such  characterization proves to be crucial in understanding the limiting behavior of the relevant estimator; see for example \cite{gjw01B}, \cite{jankowski_09} and \cite{BJRP2013}. Thus, for the sake of completeness, we give in Proposition \ref{CharLSE} below the Fenchel characterization proved in \cite{durotetal_13}, Lemma 2. Note that a typographical error occurring in \cite{durotetal_13}, Lemma 2, is now corrected. More precisely, with notation \eqref{def:F_p} and \eqref{def:H_p},
if $p \in \mathcal C$ satisfies $\sum_{k=0}^{z-1} F_{p}(k) = \sum_{k=0}^{z-1} \sum_{j=0}^k p(j) \ge \sum_{k=0}^{z-1} F_{p_n}(k)$ with equality at any knot of $p$ (instead of $\widehat p_n$ as stated in \cite{durotetal_13}) then $p = \widehat p_n$.  

\medskip

\begin{prop}\label{CharLSE}
The convex pmf $\widehat p_n$ is the LSE if and only if
\begin{eqnarray}\label{Fenchel}
H_{\hat p_n}(z)
\begin{cases}
\ge H_{p_n}(z), \ \text{for all $z \in \NN$} \\
= H_{p_n}(z), \ \text{if $z$ is a knot of $\widehat p_n$}.
\end{cases}
\end{eqnarray}
\end{prop}   

\medskip

Some remarks are in order.  The characterization above can be seen as the discrete version of the one given by \cite{gjw01B} for the LSE of a convex density with respect to Lebesgue measure. However, it should be noted that some of the consequences implied by the characterization of the continuous LSE do not hold true in our discrete case. We point here to some examples. Let $\{z_{1},\dots,z_{k}\}$ denote the set of distinct observations where we assume that $z_{1}<\dots<z_{k}$. It follows from Proposition 1 of \cite{durotetal_13} that the discrete LSE may have at most two knots in $\{z_{i}+1,\ldots,  z_{i+1}-1\} $ whenever $z_{i+1}- z_{i} > 1$.  On the other hand, Lemma 2.1 of \cite{gjw01B} shows that the continuous LSE can only bare one knot between two consecutive observations.  Another difference between the estimators stems from the lack of the notion of differentiability in the discrete case. If $F_{\hat p_n}$ and $F_{p_n}$ denote the continuous versions of the quantities defined above then $F_{\hat p_n}(s) = F_{p_n}(s)$ at any knot $s$ of the estimator; see Corollary 2.1 of \cite{gjw01B}. This equality cannot be expected to hold true for the discrete convex LSE in the general case.  In fact, by definition of $H_{p_n}$ and $H_{\hat p_n}$ in the discrete case,
$$F_{p_n}(z)=H_{p_n}(z+1)-H_{p_n}(z)  \mbox{ and } F_{\hat p_n}(z)=H_{\hat p_n}(z+1)-H_{\hat p_n}(z)$$
for all $z\in\NN$, so it follows from \eqref{Fenchel} that the equality is replaced instead by the two inequalities $F_{\hat p_n}(s) \ge F_{p_n}(s)$ and $F_{\hat p_n}(s-1) \le F_{p_n}(s-1)$.  The equality   $F_{\hat p_n}(s) = F_{p_n}(s)$ can only hold if, in addition to the equality $H_{\hat p_n}(s) = H_{p_n}(s)$, one also has $H_{\hat p_n}(s+1) = H_{p_n}(s+1)$. This happens for instance in situations where $\widehat p_{n}$ has two consecutive knots at $s$ and $s+1$.

\medskip
Next we consider almost sure consistency of $\widehat p_{n}$ in all distances $\ell_{r}$.

\medskip

\begin{prop}\label{Consis}
For any integer $r \in [2, \infty]$, we have that
\begin{eqnarray*}
\lim_{n \to \infty} \Vert \widehat p_n - p_0 \Vert_r =0
\end{eqnarray*}
with probability one.
\end{prop}

The following proposition is an easy consequence of the consistency result shown above.

\medskip

\begin{prop}\label{Knots}
If $s >0$ is a knot of $p_0$, then with probability one, there exists $n_0$ such that for all $n \ge n_0$, $s$ is a knot of $\widehat p_n$.
\end{prop}

\medskip


We finish this section by recalling boundedness in probability of $\sqrt n (\widehat p_n - p_0)$ and the implied boundedness for the associated \lq\lq integral \rq\rq \ processes. Note that boundedness in probability of $\sqrt n (\widehat p_n - p_0)$ is much weaker than tightness. However, these properties are equivalent in the case where  $\sqrt n (\widehat p_n - p_0) $ is identically equal to zero after a certain range. We will use this equivalence later on under the assumption that the true convex pmf has a finite support.  

\medskip

\begin{theorem}\label{BoundProbLSEInt}
If $\widehat p_n$ is the convex LSE of the true pmf $p_0$, then 
\begin{eqnarray}\label{BoundProbLSE}
\sqrt n \Vert \widehat p_n - p_{0} \Vert_\infty = O_p(1).
\end{eqnarray}
Furthermore, if $p_0$ has a finite support, then 
\begin{eqnarray}\label{BoundProbInt}
 \sqrt n \Vert F_{\hat p_n} - F_{p_0} \Vert_\infty = O_p(1),  \ \textrm{and} \  \sqrt n \Vert H_{\hat p_n} - H_{p_0}  \Vert_\infty = O_p(1).
\end{eqnarray}

\end{theorem}

\medskip

In the sequel, we will assume that $p_0$ has a finite support.   Finiteness of the (unknown) support of the true pmf is a reasonable assumption if we recall that the ultimate goal is estimation of species abundance.   Let $S > 0$ an integer such that $\{0, 1, \ldots, S\}$ is the support of $p_0$. Note that this means that $p_0(k) > 0$ for $k =0,1, \ldots, S$ and $p_0(k) = 0$ for $k \ge S+1$. In this sense, $S+1$ is the last knot of $p_0$ since $\Delta p_{0}(S+1)=p_0(S+2) - 2p_0(S+1) +p_0(S) = p_0(S)>0.$  Under this assumption, it is natural to ask whether the support of $\widehat p_n$ is also finite. It turns out that the answer is affirmative as we now show in the following proposition. 

\medskip

\begin{prop}\label{SuppLSE}
If $p_0$ is supported on $\{0, \ldots, S \}$ with $S \in \NN \setminus \{0\}$, then with probability one there exists $n_0$ such that for all $n \ge n_0$, the support of the LSE $\widehat p_n$ is either $\{0,\ldots,S \}$ or $ \{0, \ldots, S+1 \}$.

\end{prop}

\section{Asymptotics of the convex LSE}\label{sec:asymptotics}

In this section, we derive the weak limit of the LSE when the true distribution is supported on a finite set. The limit distribution is given in Section \ref{WeakLimit} in the most general setting where we do not make any additional assumption on the structure of the knots of the true pmf $p_0$.  It turns out that the limit distribution of the estimator involves all knots of $p_0$. The dependence on knots occurs through a stochastic piecewise convex function, $\widehat g$, that is convex between two successive knots of $p_0$; see Theorem \ref{Asymp} below.  This seems to contrast with the continuous case, where the limiting distribution of the LSE at a point  depends only on the density (and its derivatives) of the observations at this point, so that the limit distribution is ``localized''. Moreover,  the limit distribution of the MLE of a discrete log-concave pmf given in \cite{BJRP2013} is localized in some sense. For these reasons, we provide in Section \ref{Localization} below general characterizing conditions for such localizations to occur for the LSE of a discrete convex pmf. It turns out that localization does not occur in the general case but occurs, for instance, in cases where the true pmf has consecutive knots. This comprehensive study of possible localization led us to find an error in the proof of Proposition 3 in  \cite{BJRP2013}, and to conclude that the limit distribution given in Theorem 5 of that paper is not correct.

\subsection{The general setting}\label{WeakLimit}
In this section, we derive the weak limit of $\widehat p_n$ when the true distribution is supported on a finite set $\{0,\dots,S\}$. To describe the weak limit, we need to introduce the Gaussian process $\WW$ defined by
$$\WW(k) = \UU(F_{p_0}(k)) - \UU(F_{p_0}(k-1))$$
 for all integers $k=0,\dots,S+1$, where $\UU$ denotes a standard Brownian bridge. Note that $F_{p_0}(-1) =0$ which implies that $\WW(0) = \UU(p_0(0))$, and  $\WW(S+1) = 0$ since $F_{p_0}(S)=F_{p_0}(S+1)=1$. Recall that $F_{p_{n}}$ denotes the empirical distribution function corresponding to the observations $X_{1},\dots,X_{n}$. By standard results on weak convergence of empirical processes we have that
\begin{eqnarray*}
\sqrt n(F_{p_n}(k)-F_{p_0}(k)) \to_d  \UU(F_{p_0}(k))  \ \text{for $k \ \in \{0, \ldots, S+1 \}$} 
\end{eqnarray*}
where we recall that $\UU$ denotes a  standard Brownian bridge from $(0,0)$ to $(1, 0)$, and therefore,
\begin{eqnarray*}
\sqrt n(p_{n}(k)-p_{0}(k)) \to_d  \WW(k), \ \text{for $k \ \in \{0, \ldots, S+1 \}$}.  
\end{eqnarray*}
Since the LSE is described in terms of the empirical probability $p_{n}$, it is expected that the limiting distribution of the LSE can be described in terms of the limiting process $\WW$. Theorem \ref{Asymp} below proves that this is indeed the case. Precisely, the limiting distribution of the LSE is that of the minimizer (which existence is proved in Theorem \ref{Exis&Char} below) of the criterion 
\begin{eqnarray}\label{Phi}
\Phi(g)= \frac{1}{2} \sum_{k=0}^{S+1}  \left(g(k) - \WW(k)\right)^2
\end{eqnarray}
over the set $\cal C(\mathcal K)$ that we define now.  Let 
\begin{eqnarray*}
\mathcal K = \{s_1, \ldots, s_m \}  \ \textrm{or}  \ \emptyset
\end{eqnarray*}
denote the set of interior knots of $p_0$, that is all knots of $p_0$ (when they exist) that are strictly comprised between $0$ and $S+1$. For example, when $p_0$ is triangular pmf  then $\mathcal K = \emptyset$. Associated with $\mathcal K$ is the following class of functions
\begin{eqnarray*}
\mathcal C(\mathcal K) & = &  \Big \{g=(g(0), \ldots, g(S+1)) \in\RR^{S+2} \\ 
                                && \hspace{1cm} \textrm{that is convex on $\{s_j,\ldots, s_{j+1}\}$ for all  $j=0, \ldots, m$ }  \Big \}, 
\end{eqnarray*}
with  $s_{0}:=0$ and $s_{m+1}:=S+1$ (the greatest knot of $p_{0}$). This means that any $g\in\mathcal C(\mathcal K)$ is convex at points $1,\dots,S$ with possible exception at points in $\mathcal K$.

\medskip

Using the notation \eqref{eq:defH}, we have the following theorem. 

\medskip

\begin{theorem}\label{Exis&Char}
The criterion (\ref{Phi}) admits a unique minimizer $\widehat g$ over $\cal C(\cal K)$. Furthermore, an element $\widehat g\in\cal C(\cal K)$ is the minimizer if and only if the process $ \widehat{\HH}$ defined on $\{0,\ldots, S+2\}$ by 
\begin{eqnarray}\label{hatH}
\widehat\HH(x) = 
\begin{cases}
\sum_{k=0}^{x-1} \sum_{j=0}^k \widehat g(j), \ \text{if $x \in \{1, \ldots, S+2\}$} \\
\HH(0)= 0, \ \hspace{1.25cm} \text{if $x=0$}
\end{cases}
\end{eqnarray}
satisfies
\begin{eqnarray}\label{CharLSIneq2}
\widehat \HH(x)
\begin{cases}
\ge  \HH(x), & \ \textrm{for all} \ x \in \{0, \ldots, S+2\}  \\
=  \HH(x), & \ \textrm{if  $x \in \{s_{0},\dots,s_{m+1}, S+2\}$ }  \\
                 &  \hspace{0.2cm} \textrm{or $x$ is a knot of $g$ in $\{s_{j}+1,\dots,s_{j+1}-1\}$} \\
                & \hspace{0.2cm} \textrm{for some $j=0,\dots,m$}.
\end{cases}
\end{eqnarray}

\end{theorem}

\medskip

In the above theorem, note that it is implicit that the minimizer $\widehat g$ is actually $\widehat g(\omega)$. We are now ready to state the result of weak convergence of $\widehat p_n$. For $x \in \{0, \ldots, S+1 \}$ define
\begin{eqnarray*}
\widehat{\mathbb G}(x) = \sum_{k=0}^{x} \widehat g(k) = \widehat \HH(x+1) - \widehat \HH(x).
\end{eqnarray*}
Recall that with probability one, $\hat p_{n}(j)=p_{0}(j)=0$ for all $j>S+1$ provided that $n$ is sufficiently large, see Proposition \ref{SuppLSE}. Therefore, it suffices to compute the weak limit  of $\sqrt n(\hat p_{n}-p_{0})$ on $\{0,\dots,S+1\}$.

\medskip

\begin{theorem}\label{Asymp}
If $\widehat p_n$ is the convex LSE of the true pmf $p_0$ with support $\{0, \ldots, S \}$, then we have the joint weak convergence on $\{0,\dots,S+1\}$
\begin{eqnarray*}
\left (
\begin{array}{lll}
\sqrt n (H_{\hat p_n} - H_{p_0}) \\
\sqrt n (F_{\hat p_n} - F_{p_0})  \\
\sqrt{n} (\widehat p_n - p_0)   
\end{array}
\right) \Rightarrow
\left ( 
\begin{array}{lll}
\widehat \HH \\
\widehat{\mathbb G} \\
\widehat g
\end{array}
\right),
\end{eqnarray*}
as $n \to \infty$. 
\end{theorem}
 
\subsection{Localization}\label{Localization}
In the previous subsection, it has been proved that the limiting distribution of $\widehat p_{n}$ at a fixed point involves all knots of $p_{0}$ in the general case; see Theorem  \ref{Asymp}. This seems to contrast with the continuous case where the  limiting distribution of the LSE is localized in the sense that it depends only on the true density (and its derivatives) at the fixed point. 
A natural question that arises from the above consideration, and which we answer below, is on whether the convex LSE could be localized in the discrete case. 

\medskip

To draw a correct comparison between what happens in the discrete and continuous cases, one has to go back to the working assumptions under which the limiting distribution has been derived in the latter case. In \cite{gjw01B} it is assumed that the true convex density $f$ defined on $[0, \infty)$ is twice continuously differentiable in a small neighborhood of a fixed point $x_0 > 0$ such that $f''(x_0) > 0$. In particular, the density is strictly convex at $x_{0}$. Under this assumption,  the limiting distribution  is $\sigma(x_0)^{-1} Y$ where $\sigma(x_0)  = 24^{1/5} [f(x_0)]^{-2/5} [f''(x_0)]^{-1/5}$  and $Y$ has the same distribution as the second derivative at zero of the outer envelope of the first integral of a two-sided Brownian motion plus the drift $t \mapsto t^4$. See \cite{gjw01A} for a complete analysis of existence and uniqueness of this outer envelope or \lq \lq invelope\rq\rq \ as has been coined by the authors. In our Theorem \ref{Asymp} we do not consider any particular configuration for the knots of $p_0$. For the sake of comparison, and if we translate for the moment strict convexity at a point $s$ in the support of $p_0$ as having $s$ to be \textrm{a triple knot}, that is $s-1$, $s$ and $s+1$ are successive knots of $p_0$, then it follows from Proposition \ref{Knots} and Proposition \ref{CharLSE} that with probability one there exists $n_0$ large enough such that for all $n \ge n_0$  
\begin{eqnarray*}
H_{\hat p_n}(s-1) = H_{p_n}(s-1), \ \ H_{\hat p_n}(s) = H_{p_n}(s), \ \ \textrm{and} \ \  H_{\hat p_n}(s+1) = H_{p_n}(s+1). 
\end{eqnarray*}
This implies that 
\begin{eqnarray}\label{TripleKnot}
\widehat p_n(s) = p_n(s)
\end{eqnarray}
and the limiting distribution is simply that of the Gaussian random variable $\mathcal {N}(0, p_0(s)(1 - p_0(s))$. In this case, the limit of the LSE at the point $s$ is completely localized in the sense that it is not influenced by the remaining knots of $p_0$. The identity in (\ref{TripleKnot}) shows even the stronger fact that the localization is actually happening at the level of the estimator itself. 

\medskip 
It is conceivable that other configurations lead to some form of localization of the weak limit. We provide below general characterizing conditions for such localizations to occur. In fact,  the LSE and its weak limit get localized either to the left or right at any knot of $p_0$ that is either followed or preceded by another knot of $p_0$. In such cases, the limit of the LSE can be described only in terms of knots of $p_0$ that are either to the left or right of that knot; see  comments after Theorem \ref{LocalLeft} and Theorem \ref{LocalRight}. In the sequel, we shall use the same notation as in Section \ref{WeakLimit}.

\medskip

First, we consider the question of localizing ``to the left'' of a knot. This means that given a knot $s\in\{s_{1},\dots,s_{m}\}$, we wonder whether the restriction to $\{0,\dots ,s\}$ of the limiting $\widehat g$ is distributed as the minimiser of the left-localized criterion
$$\Phi^{\leq s}(g) = \frac{1}{2} \sum_{k =0}^s (g(k)- \WW(k))^2$$
over the set 
\begin{eqnarray*}
{\mathcal C}^{\leq s}(\mathcal K) & = &  \Big \{g=(g(0), \ldots, g(s)) \in \RR^{s+1} \\ 
                                && \hspace{1cm} \textrm{that is convex on $\{s_j,\ldots, s_{j+1}\}$} \\
                                && \hspace{1cm} \textrm{for all $j \in \{0, \ldots, m\}$ with $s_{j+1}\le s$ } \Big\}.
\end{eqnarray*}
The following theorem provides a necessary and sufficient condition for the answer to be positive. It also gives a necessary and sufficient condition for the restriction of $\sqrt n(\widehat p_{n}-p_{0})$ to $(0,\dots,s)$ to converge to the left-localized minimizer.

\medskip

\begin{theorem}\label{LocalLeft}
Assume that the support of $p_{0}$ is finite. Then for an arbitrary $s\in\{s_{1},\dots,s_{m}\}$, there exists a unique minimizer of $\Phi^{\leq s}$ over ${\mathcal C}^{\leq s}(\mathcal K)$. The minimizer is equal to $(\widehat g(0),\dots,\widehat g(s))$ if, and only if, 
\begin{equation}\label{NSC Left}
\widehat{\mathbb G}(s)=  \UU(F_{p_0}(s)).
\end{equation}
Moreover, $\sqrt n(\widehat p_{n}(0)-p_{0}(0),\dots,\widehat p_{n}(s)-p_{0}(s))$ converges in distribution to the minimizer of $\Phi^{\leq s}$ over ${\mathcal C}^{\leq s}(\mathcal K)$ if, and only if,
\begin{equation}\label{NSC Left n}
F_{\hat p_n}(s)=F_{p_n}(s)+o_{p}(n^{-1/2}).
\end{equation}
\end{theorem}
It is worth mentioning that \eqref{NSC Left n} holds for instance if $s$ is a double knot, in the sense that both $s$ and $s+1$ are knots of $p_{0}$. Indeed, it follows from Proposition \ref{Knots} together with the characterization in Proposition \ref{CharLSE}, that if $s$ is a double knot of $p_{0}$, then with probability one, both $H_{\hat p_n}(s)=H_{p_n}(s)$ and $H_{\hat p_n}(s+1)=H_{p_n}(s+1)$ hold true for sufficiently large $n$. Therefore, 
$$F_{\hat p_n}(s)=F_{p_n}(s)$$
with probability one for sufficiently large $n$, so that \eqref{NSC Left n} holds and the limiting distribution is left-localized. 

\medskip


\medskip

Now, we consider the question of localizing ``to the right'' of a knot. This means that given a knot $s\in\{s_{1},\dots,s_{m}\}$, we wonder whether the restriction to $\{s,\dots ,S+1\}$ of the limiting $\widehat g$ is distributed as the minimizer of the right-localized criterion
$$\Phi^{\geq s}(g) = \frac{1}{2} \sum_{k =s}^{S+1} (g(k)- \WW(k))^2$$
over the set 
\begin{eqnarray*}
{\mathcal C}^{\geq s}(\mathcal K) & = &  \Big \{g=(g(s), \ldots, g(S+1)) \in \RR^{S-s}\\ 
                                && \hspace{1cm}  \textrm{ that is convex on $\{s_j,\ldots, s_{j+1}\}$} \\
                                && \hspace{1cm} \textrm{for all $j \in \{0, \ldots, m\}$ with $s_{j}\geq s$}  \Big \}.
\end{eqnarray*}
A necessary and sufficient condition for the answer to be positive, is given below.

\medskip

\begin{theorem}\label{LocalRight}
Assume that the support of $p_{0}$ is finite. Then for an arbitrary $s\in\{s_{1},\dots,s_{m}\}$, there exists a unique minimizer of $\Phi^{\geq s}$ over ${\mathcal C}^{\geq s}(\mathcal K)$. The minimizer is equal to $(\widehat g(s),\dots,\widehat g(S+1))$ if, and only if, 
\begin{equation}\label{NSC Right}
\widehat{\mathbb G}(s-1)=  \UU(F_{p_0}(s-1)).
\end{equation}
Moreover, $\sqrt n(\widehat p_{n}(s)-p_{0}(s),\dots,\widehat p_{n}(S+1)-p_{0}(S+1))$ converges in distribution to the minimizer of $\Phi^{\geq s}$ over ${\mathcal C}^{\geq s}(\mathcal K)$ if, and only if,
\begin{equation}\label{NSC Right n}
F_{\hat p_n}(s-1)=F_{p_n}(s-1)+o_{p}(n^{-1/2}).
\end{equation}
\end{theorem}
Similar as above,  \eqref{NSC Right n} holds in the specific case where  $s-1$ is a double knot in the sense that both $s-1$ and $s$ are knots of $p_{0}$. 
\medskip

\medskip

Theorems \ref{LocalLeft} and \ref{LocalRight} above give some insight on possible left or right localization of the limiting distribution. For the purpose of implementing the estimator $\widehat p_{n}$, it is of interest not only to localize the limiting distribution, but also to localize the estimator itself. To be more specific, suppose that we are not interested in estimating the whole pmf $p_{0}$, but that we are interested only in the estimation of $p_{0}(k_{1}),\dots,p_{0}(k_{I})$ for given $k_{1}<\dots<k_{I}$. Of course, one can nevertheless compute the whole estimator $\widehat p_{n}$ as the minimizer of the criterion $\Phi_{n}$ over the set of convex sequences with a finite $\ell_{2}$-norm, and then extract the estimators $\widehat p_{n}(k_{1}),\dots,\widehat p_{n}(k_{I})$, but it may be more convenient to minimize instead a localized criterion. The following theorem give conditions for the localization to be allowed: it proves that in $\Phi_{n}$, the sum over $\NN$ can be replaced by a finite sum under appropriate conditions. The proof involves similar arguments as the proofs of Theorems \ref{LocalLeft} and \ref{LocalRight} combined to Proposition \ref{Knots},  so it is omitted.

\begin{theorem}\label{LocalBoth}
Assume that $p_{0}$ is supported on $\{0,\dots,S\}$.
Let $s \leq k_{1}$ be such that either $s-1$ is a double knot (if exists) of $p_{0}$, or $s=0$. Let $s'\geq k_{I}$ be either a double knot (if exists) of $p_{0}$, or $s'\geq S+1$. Let $z\leq s$ and $z'\geq s'$ be arbitrary points. Then, with probability one, there exists $n_{0}$ such that for all $n\geq n_{0}$,
$$\Big(\widehat p_{n}(k_{1}),\dots,\widehat p_{n}(k_{I})\Big)$$
is the restriction to $\{k_{1},\dots,k_{I}\}$ of the minimizer of 
\begin{equation}\label{CritLocal}
\frac12\sum_{k=z}^{z'}(p_{n}(k)-p(k))^2
\end{equation}
over the set of convex sequences $(p(z),p(z+1),\dots,p(z')).$ Moreover, 
$$\sqrt n\Big(\widehat p_{n}(k_{1})-p_{0}(k_{1}),\dots,\widehat p_{n}(k_{I})-p_{0}(k_{I})\Big)$$
converges in distribution to the restriction to $\{k_{1},\dots,k_{I}\}$ of the minimizer of 
$$\frac12\sum_{k=z}^{z'}(g(k)-\WW(k))^2$$
over the set of sequences $g=(g(z),\dots,g(z'))$ that are convex on $\{s_j,\dots, s_{j+1}\}$   for all $j=0, \ldots, m$ with $s_{j}\geq z$ and $s_{j+1}\leq z'$.
\end{theorem}
This means that instead of computing $\hat p_{n}(k_{1}),\dots,\hat p_{n}(k_{I})$, one can compute the minimizer of \eqref{CritLocal}
over the set of convex sequences $(p(z),p(z+1),\dots,p(z')),$ and then take the restriction to $\{k_{1},\dots,k_{I}\}:$ the restriction has the same limit distribution as $\hat p_{n}(k_{1}),\dots,\hat p_{n}(k_{I})$, and the two vectors are even equal with probability one for sufficiently large $n$.

Note that  the precise location of the knots, double knots, and $S$, need not to be known to compute $\widehat p_{n}(k_{1}),\dots,\widehat p_{n}(k_{I})$ using the localized criterion \eqref {CritLocal}. However, $z$ and $z'$ have to be chosen small enough, and large enough respectively, to ensure that the localization is allowed.

\section{Numerical aspects}\label{sec:implementation}

\subsection{A Dysktra algorithm for computing the asymptotic distribution}\label{DykstraAlgo}

In monotone or (concave/convex) regression, active set methods are often proposed to compute the estimators; see for example~\cite{HF85} or~\cite{GJW08} for more recent work. On may also refer to~\cite{MAM91} for a more general method of non parametric regression under shape constraint.
 
Here, we describe a simple algorithm that enables us to simulate a sample of any size from the asymptotic distributions of Theorem \ref{Asymp}. Here, we focus on the weak limit of the convex LSE itself. As shown above, $\sqrt n (\widehat p_n - p_0) \Rightarrow \widehat g$. When the setting is that of a convex pmf with a finite support $\{0, \ldots, S\}$, it follows from Theorem \ref{Exis&Char} above that 
\begin{eqnarray*}
\widehat g = \text{argmin}_{g \in \mathcal{C}(\mathcal K)}  \Phi(g) = \text{argmin}_{g \in \mathcal{C}(\mathcal K)}  \ \frac{1}{2} \sum_{k=0}^{S+1} \left( g(k) - \WW(k)\right)^2.
\end{eqnarray*}
where $\mathcal{K} = \{s_1, \ldots, s_m\}$ is the set of interior knots of $p_0$ and $\mathcal C(\mathcal K)$ the class of piecewise functions $g$ defined on $\{0, \ldots, S+1\}$ such that $g$ is convex on $\{s_j, \ldots, s_{j+1} \}$ for $j=0, \ldots, m$ with $s_0\equiv 0$ and $s_{m+1} \equiv S+1$. Hence, if we denote $\mathcal{K}_j = \{s_j, \ldots, s_{j+1} \}$, then  $\widehat{g}$ can then be viewed as the projection of  $\WW$ onto intersection of the convex subsets
\begin{eqnarray*}
\mathcal{C}_{j} = \left\{ g \ \textrm{is piecewise on $\{0, \ldots, S+1\}$}: \  g \mbox{ is convex on $\mathcal{K}_j$} \right\}
\end{eqnarray*}
for $j=0, \ldots, m$. It is easily seen that $\mathcal{C}_{j}, j =0, \ldots, m$ are closed convex cones. Therefore, minimizing $\Phi$ over
$\mathcal C(\mathcal K)$ is equivalent to minimizing it over $\bigcap_{j=0}^m \mathcal{C}_j.$  The solution can be found using the algorithm of  \cite{dykstra1983algorithm} which proceeds by performing cyclic projections onto the convex cones $\mathcal{C}_0, \ldots, \mathcal{C}_m$.
Although details of the algorithm is given in \cite{dykstra1983algorithm}, we describe here how these projections are performed.  At iteration $0$, set 
\begin{eqnarray*}
g^{(0)}  \equiv (\WW(0), \WW(1), \ldots, \WW(S+1)), \ \textrm{and the increments} \  u^{(0)}_j = 0 \ \ \textrm{for $j= 0,...,m$}.
\end{eqnarray*}
For $ n \ge 1$, the algorithm will then proceed in the three following steps.
 
\medskip

\begin{description}
\item [(1)]  Compute $g_{j}^{(n)}$, the projection of $g_{j-1}^{(n)} - u_{j}^{(n-1)}$  onto $\mathcal{C}_{j}$ for $j =0, \ldots, m$.

\item [(2)]  Set $u_{j}^{(n)} \equiv g_{j}^{(n)} - \left(g_{j-1}^{(n)} - u_{j}^{(n-1)}\right)$.

\item [(3)] Set $n = n+1$ and go to $\textbf{(1)}.$

\end{description}

\medskip

Granted that we know how to obtain the convex projections $g^{(n)}_j$, convergence of the above algorithm is a consequence of Theorem $3.1$ of ~\cite{dykstra1983algorithm}. The projection onto each cone $\mathcal{C}_j, j =0, \ldots, m$ can be efficiently computed using the \proglang{R} function  {\tt conreg} available in the \proglang{R} package \pkg{COBS}; see \cite{NM07} for more details.

\medskip

Now, for any fixed integer $N \ge 1$, a sample of size $N$ from the same distribution as $(\widehat g(0), \ldots, \widehat g(S+1))$ can be done as follows: we generate a centered Gaussian vector $(W_0, \ldots, W_S)$ with dispersion matrix $\Sigma$ given by 
$$\Sigma_{ij} = p_0(i)  \mathbb{I}_{i = j} - p_0(i) p_0(j)$$ for $ 0 \le i, j \le S$. This can be done using the  \proglang{R} function {\tt rmvnorm} available from the \pkg{mvtnorm} package. In the second step, we compute the piecewise convex projection of $(W_0, \ldots, W_S, 0)$ as described above, and the two steps are then repeated $N$ times.



\subsection{How well the true knots are captured}\label{CaptureKnots}
Recall that Proposition \ref{Knots} implies that with probability one, having enough large sample sizes ensures that a knot of the true pmf $p_0$ is also a knot of the LSE $\widehat p_n$. However, the proposition does not indicate how much large $n$ should be. To gain some insight into the relationship between the size of the sample at hand and whether the knots of the estimator include all true knots, we have carried out a simulation study with samples of size $n \in \{50, 200, 800, 3200, 12800, 51200\}$.  Given a simulated sample of size $n$ from a distribution $p_{0}$,  the convex LSE $\hat p_{n}$ was computed using the algorithm described in \cite{durotetal_13}.

\medskip

To define the convex pmf's under which the samples were generated, we need a few notation. 
Given $j \in \NN \setminus \{0\}$, consider the triangular pmf with support on $\{0, \ldots, j-1 \}$ given by 
$$
T_j(i)   = \frac{2(j-i)_+}{j(j+1)}  .
$$
As proved in \cite{durotetal_13}, Theorem 7,  a pmf $p_{0}$ is convex if and only if $p_{0}$ admits the mixture representation
\begin{eqnarray*}
p_{0} = \sum_{j \ge 1} \pi_j T_j
\end{eqnarray*}
where $0 \le \pi_j \le 1$ and $\sum_{j\ge 1} \pi_j =1$. The representation is unique and the mixing weights are given by 
\begin{eqnarray*}
\pi_j = \frac{j(j+1)}{2} \Delta p(j)  = \frac{j(j+1)}{2} (p_{0}(j+1)+ p_{0}(j-1) - 2 p_{0}(j)) 
\end{eqnarray*}
for $j \ge 1$.  Note that it follows from the expression of $\pi_j$ that $j$ is a knot of $p_{0}$ if and only if $\pi_j > 0$. In particular, if the support of $p_{0}$ takes the form $\{0,\dots,S\}$, then $\pi_{S+1}>0$ whereas $\pi_{j}=0$ for all $j> S+1$.

\medskip

In our simulations, the samples were generated from four convex pmf's that are all supported on $\{0, \ldots, 10 \}$. We give in Table \ref{MixingWeights} the values of the mixing weights $\pi_j, \ 1 \le j \le 11$ for various pmf's $p_{0}$ that are denoted by $p_{1}$, $p_{2}$, $p_{3}$ and $p_{4}$. 
\medskip

\begin{table}[!h]
\caption{Mixing weights $\pi_j$ for the convex pmf's $p_1, p_2, p_3, p_4$. } 
\vspace{0.4cm}
\centering 
\begin{tabular}{|c||c|c|c|c|c|c|c|c|c|c|c|}
\hline\hline 
pmf & $\pi_1$ & $\pi_2$ & $\pi_3$  & $\pi_4$ & $\pi_5$  & $\pi_6$ & $\pi_7$ & $\pi_8$ & $\pi_9$ & $\pi_{10}$ & $\pi_{11}$ \\ 
\hline
\hline
$p_1$ &   0  & 1/6  & 0  & 0&  1/6 &  0 & 0 & 0  & 1/2 & 0 & 1/6 \\ 
\hline
$p_2$ &   0  & 0  & 0  & 1/6 &  0 &  1/6 & 0 & 1/12  & 0 & 1/2 & 1/12 \\ 
\hline
$p_3$ &   0  &0  & 1/6  & 1/12 &  1/4 &  0 & 1/12 & 0  & 1/6 & 1/6 & 1/6 \\ 
\hline
$p_4$ &   0  & 1/12  & 1/6  & 1/12 &  1/12 &  1/12 & 1/12 & 1/12  & 1/12 & 1/6 & 1/12 \\ 
\hline
\hline
\end{tabular}
\label{MixingWeights}
\end{table}

\begin{figure}[h!]
\begin{center}
\includegraphics[clip=true, width=7.3cm]{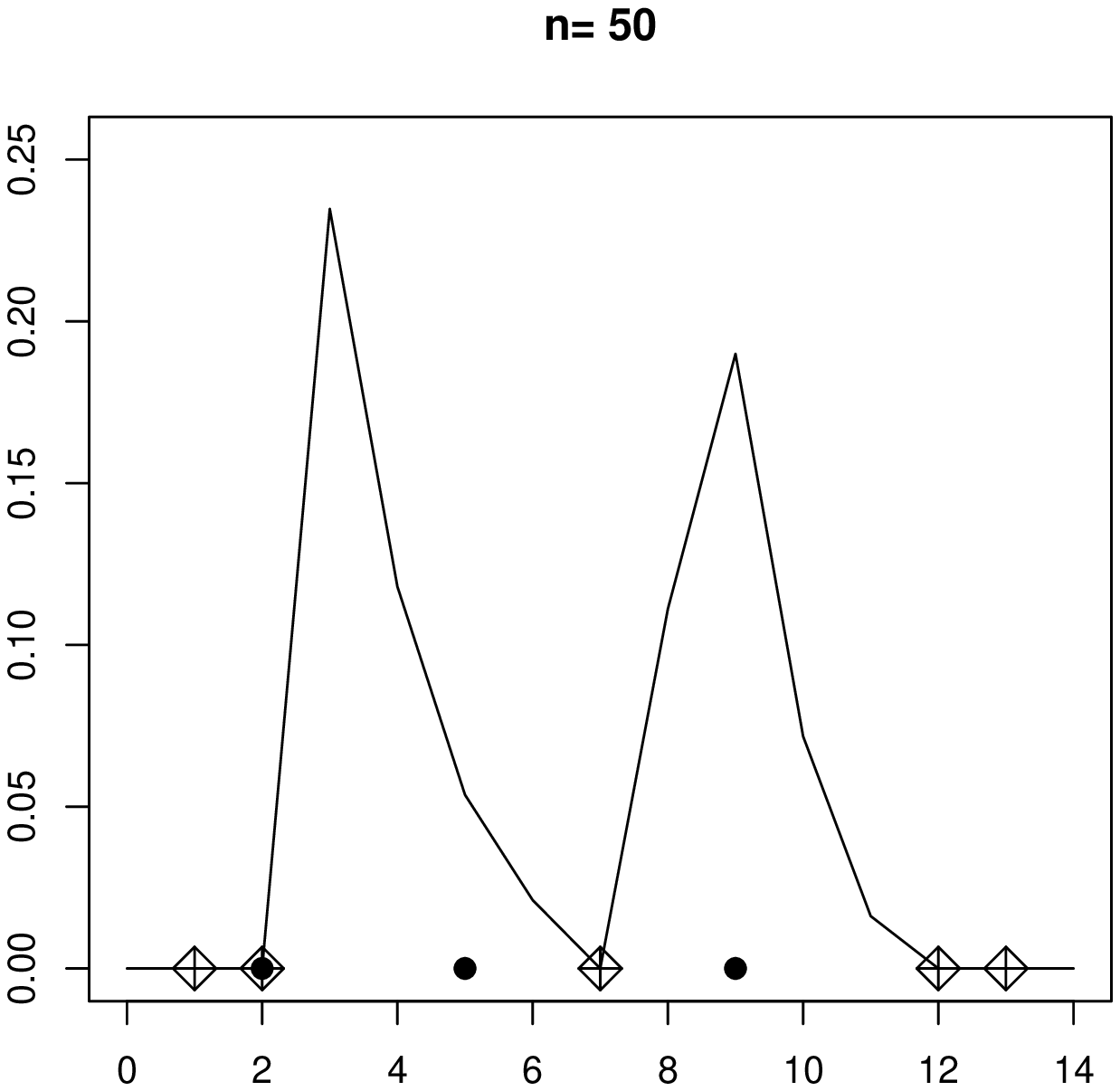}  \includegraphics[clip=true, width=7.3cm]{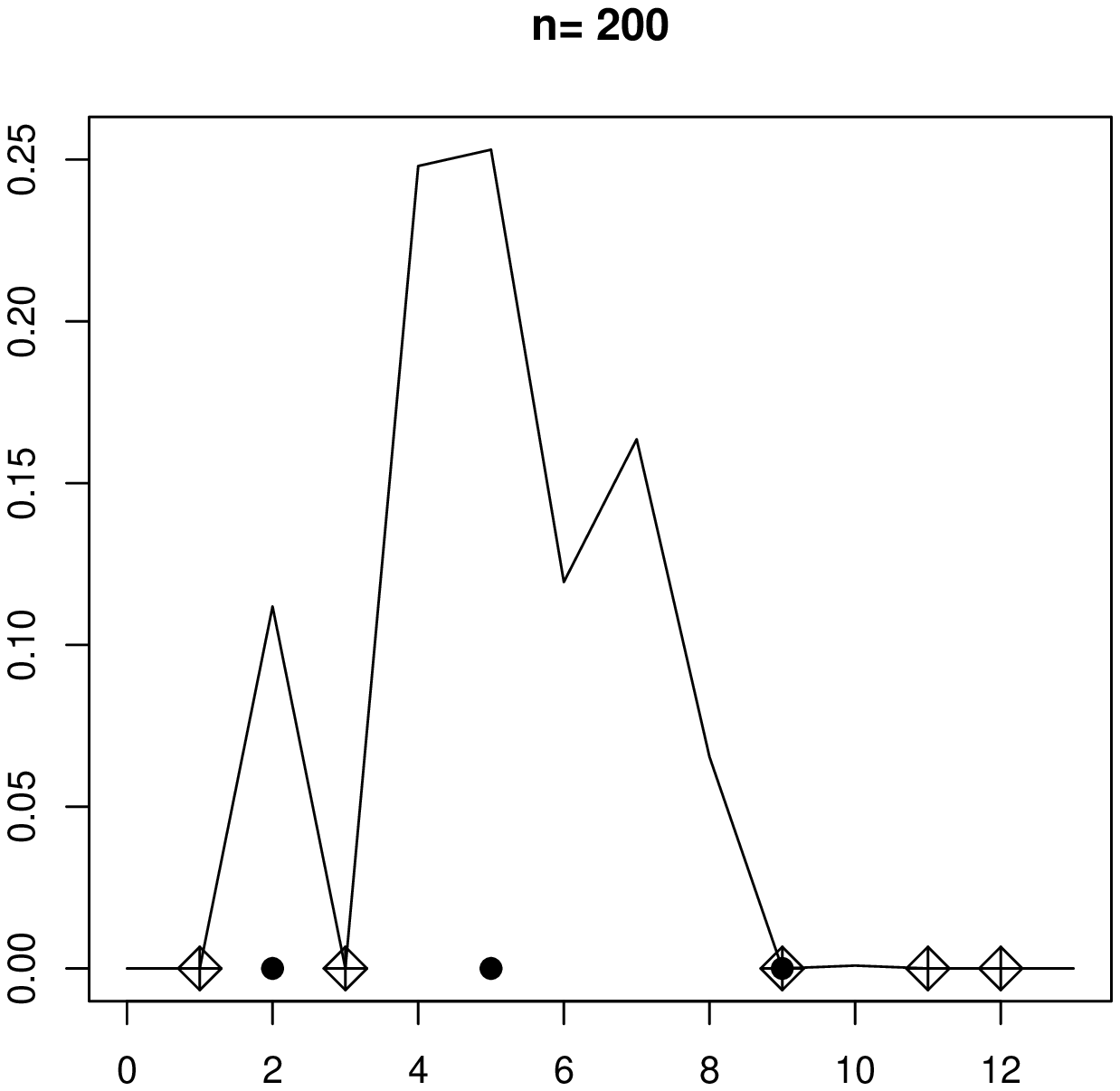}\\
\vspace{-1.3cm}
\end{center}
\caption{The figures show the process $\sqrt n (H_{\hat p_n} - H_{p_n})$ for a random sample of size $n$ as shown. The samples were generated from a convex pmf supported on $\{0,1, \ldots, 10 \}$ with interior knots at $2, 5, 9$ and mixing probabilities $\pi_2=\pi_5=\pi_{11}=1/6$, $\pi_9=1/2$ and $\pi_j = 0$ for $j \in \{1,3,4,6,7,8\}$ (see text for details). The diamond symbols depict the knots of the LSE $\widehat p_n$ computed based on the samples, whereas the bullets show the locations of the true knots.}
\label{CharT2}
\end{figure}

\begin{figure}[h!]
\begin{center}
\includegraphics[clip=true, width=7.3cm]{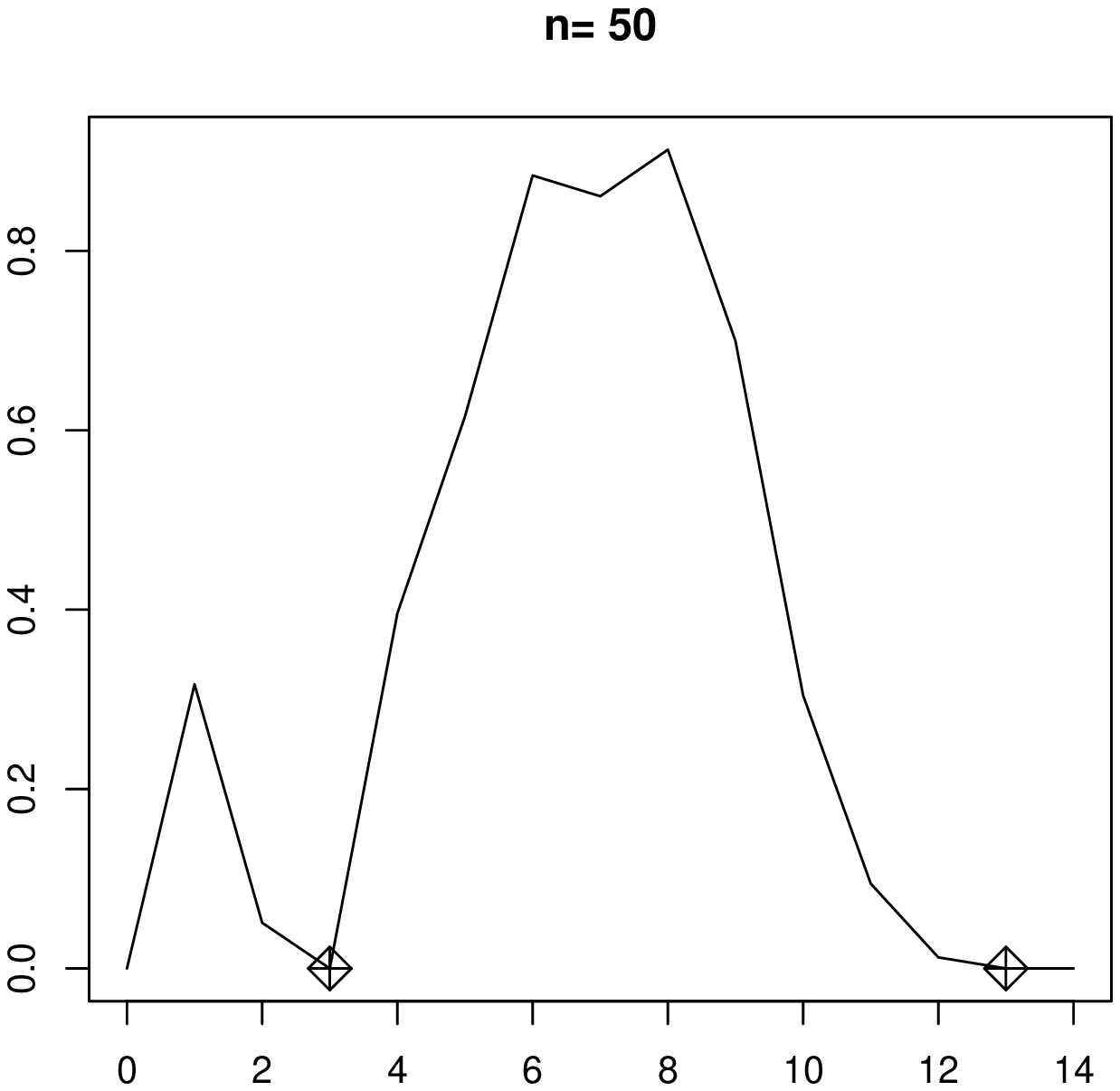}  \includegraphics[clip=true, width=7.3cm]{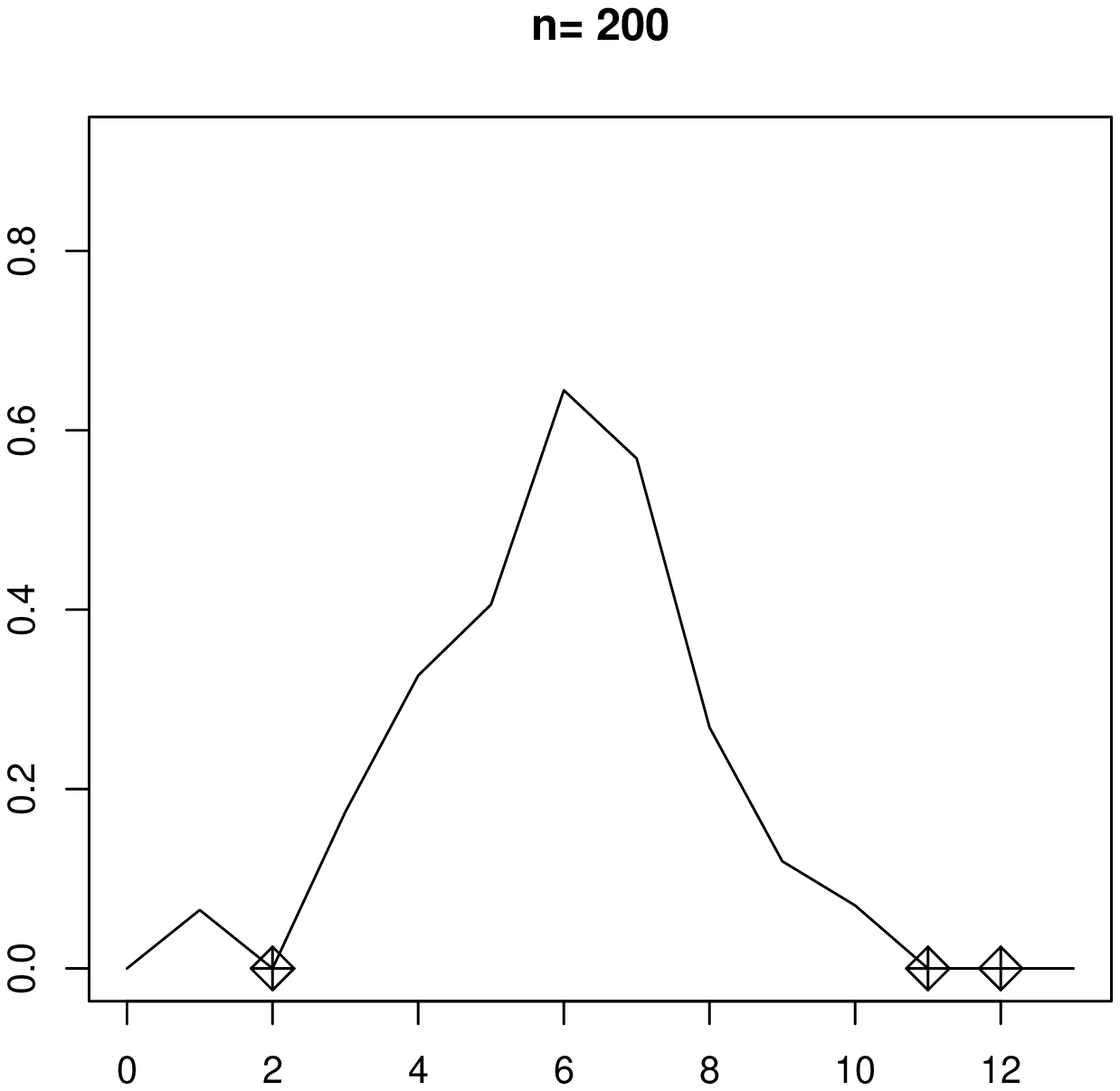}\\
\vspace{-1.3cm}
\end{center}
\caption{The figures show the process $\sqrt n (H_{\hat p_n} - H_{p_n})$ for a random sample of size $n$ as shown. The samples were generated from a triangular pmf supported on $\{0,1, \ldots, 10 \}$. The diamond symbols depict the knots of the LSE $\widehat p_n$ computed based on the samples.}
\label{CharT}
\end{figure}

\medskip
Figure \ref{CharT2} shows the process $\sqrt n (H_{\hat p_n} - H_{p_n})$ together with the knots of the LSE $\hat p_{n}$ and the true knots, for a sample of size $n\in\{50,200\}$ generated from $p_{1}$. On these examples, it can be seen that, in accordance with Proposition \ref{CharLSE}, $H_{\hat p_n}(z) \geq H_{p_n}(z)$ with an equality at all knots $z$ of $\hat p_{n}$. However, the sample sizes are not large enough to ensure that the knots of $\hat p_{n}$ include all knots of the true pmf $p_{1}$. Neither they are large enough to ensure that the support of $\hat p_{n}$ is included in $\{0,\dots,S+1\}$ where $S=10$ denotes the greatest point in the support of the true pmf; see Proposition \ref{SuppLSE}. Figure \ref{CharT} is similar to Figure \ref{CharT2} but now, the samples are generated from the triangular pmf $T_{11}$, which means that the mixing probabilities are $\pi_{11}=1$ and $\pi_{j}=0$ for all $j\neq 11$. Again, we observe that $H_{\hat p_n}(z) \geq H_{p_n}(z)$ with an equality at all knots $z$ of $\hat p_{n}$. In the case of a sample size $n=50$, the knots of $\hat p_{n}$ do not include the only true knot 11 and the support of $\hat p_{n}$ is not included in $\{0,\dots,11\}$.  On the other hand, in the case of a larger sample size $n=200$, the knots of $\hat p_{n}$ include the only true knot 11 and $\hat p_{n}$ is supported on $\{0,\dots,11\}$.

Figures \ref{CharT2} and \ref{CharT} are not sufficient to gain some insight into the connection between the knots of $\hat p_{n}$ and the true knots since only one sample is considered in each situation. Thus, for each considered sample size and distribution, we simulated independently 1000 samples to evaluate the probability that the knots of $\hat p_{n}$ include all true knots. The probabilities are estimated by empirical frequencies. Results are reported in Table \ref{table}. As expected, the empirical frequency increases as $n$ increases. It is typically larger in cases of true distributions with only few knots than in cases of true distributions with many knots.

\begin{table}[!h]
\caption{Empirical frequencies in $\%$ of having \textrm{all} knots of the true convex pmf among those of the estimator $\widehat p_n$ for $n \in \{50, 200, 800, 3200, 12800, 51200\}$. The empirical frequencies are based on 1000 replications for each sample size and distribution. The true convex pmf's, $p_1$, $p_2$, $p_3$ and $p_4$, have $3$, $4$, $6$ and $9$ interior knots respectively. See text for the exact expressions of those pmf's.} 
\vspace{0.4cm}
\centering 
\begin{tabular}{|c||c|c|c|c|} 
\hline\hline 
$n$ & $p_1$ & $p_2$ & $p_3$  & $p_4$ \\ 
\hline
\hline 
50 &   4.2 & 0.0  &  0.0 & 0.0   \\ 
\hline 
200 &   14.5 &  2.1 & 0.3  &   0.0    \\ 
\hline
800 &   47.3  &  9.0 &  3.7 &   0.0    \\ 
\hline
3200 &   84.1 &  31.4 &  32.4 &  4.3    \\ 
\hline
12800 &   99.1 &  66.3 &  71.2 &  33.1     \\ 
\hline
51200 & 100   & 92.8  & 95.9 &     88.5  \\ 
\hline
\hline
\end{tabular}
\label{table}
\end{table}

\subsection{Assessing the convergence of the convex LSE to the weak limit}\label{AssessConv}

To assess convergence of the estimation error to the right weak limit, consider $\widehat{\mathbb{F}}^{(j)}_{n,M}$ and $\mathbb {F}^{(j)}_{M'}$ to be respectively the empirical distributions of $\sqrt n (\widehat p_n(j) - p_0(j))$ and $\widehat g(j)$  for $j \in \{0, \ldots, S+1\}$ based on $M$ and $M'$ independent replications. Here, $M$  and $M'$ will be chosen to be large. More explicitly, a sample of  $n$ independent random variables $X^{(i)}_1, \ldots, X^{(i)}_n$  is drawn from $p_0$ for each $ i =0, \ldots, M$ to form a sample of size $M$ from the distribution of the estimation error.  Note that this sample is multidimensional of dimension $S+1$ hence our need to consider the marginal components of its distribution. Similarly, we draw a sample of size $M'$ from the distribution of weak limit using  the algorithm described in Section \ref{DykstraAlgo}.  Define now
\begin{eqnarray*}
D_{n,M, M'} = \sup_{0 \le j \le S+1}  \big \Vert \widehat{\mathbb{F}}^{(j)}_{n,M}  - \mathbb F^{(j)}_{M'} \big \Vert_\infty.
\end{eqnarray*}
We will use this random variable to assess the established convergence. This is based on the fact that it is expected to become small for large $n$. Since the \lq\lq target \rq\rq \ distributions are  $\mathbb F^{(j)}_{M'}, j=0, \ldots, S+1$, we choose $M' > M$. Also, we store those obtained empirical distributions and re-use them while sampling many times from the estimation error. This enables us to obtain independent realizations from $\widehat{\mathbb{F}}^{(j)}_{n,M}$ while $\mathbb F^{(j)}_{M'}, j =0, \ldots, S+1$ are fixed.  To visualize the statistical summary of $D_{n,M, M'}$, the obtained outcomes are represented in the form of boxplots. Those were based on 100 replications of  $\widehat{\mathbb{F}}^{(j)}_{n,M}, j =0, \ldots, S+1$. Here, $M=1000$, $M'=5000$ and $n \in \{50, 100, 500, 1000, 5000, 10000, 500000 \}$.  In the simulations, we have taken the following true convex pmf's which are all supported on $ \{0, \ldots, 10 \}$:

\begin{itemize}

\item The triangular pmf  $p_0$, that is 
\begin{eqnarray*}
p_0(i) = \left \{
\begin{array}{ll}
\frac{(11-i)}{66}, \ i \in \{0, \ldots, 10 \} \\
0, \hspace{0.9cm} \textrm{otherwise}.
\end{array}
\right.
\end{eqnarray*}

\item The convex pmf's $p_1, p_2, p_3, p_4$ considered above in Section \ref{CaptureKnots}. See also Table \ref{MixingWeights} above. 

\item The pmf, $p_5$, of a truncated Geometric pmf with success probability equal to $1/2$, that is 
\begin{eqnarray*}
p_5(i) =\left\{
\begin{array}{ll}
 (1- 2^{-11})^{-1}\  2^{-(i+1)},  \ \  i \in \{0, \ldots, 10 \} \\
0, \hspace{3.6cm} \textrm{otherwise}.
\end{array}
\right.
\end{eqnarray*}

\end{itemize}

Note that if a geometric pmf is always convex on $\NN$, this is not the case anymore after truncation. Indeed, convexity of the latter version holds true if and only if the waiting probability is $\le 1/2$. A simple proof of this fact can be found in Section \ref{sec:proofs}. 

\medskip

\begin{lemma}\label{lem:geometric}
Let $S$ be a positive integer and $q\in(0,1)$. The truncated Geometric distribution, defined by
$$p(i)=\frac{q^i(1-q)}{1-q^{S+1}},\ i\in\{0,\dots,S\}$$
and $p(i)=0$ for all integers $i\geq S+1$, is convex if and anly if $q\le 1/2$.
\end{lemma}

\medskip

\begin{figure}[h!]
\begin{center}
\includegraphics[clip=true, width=9cm]{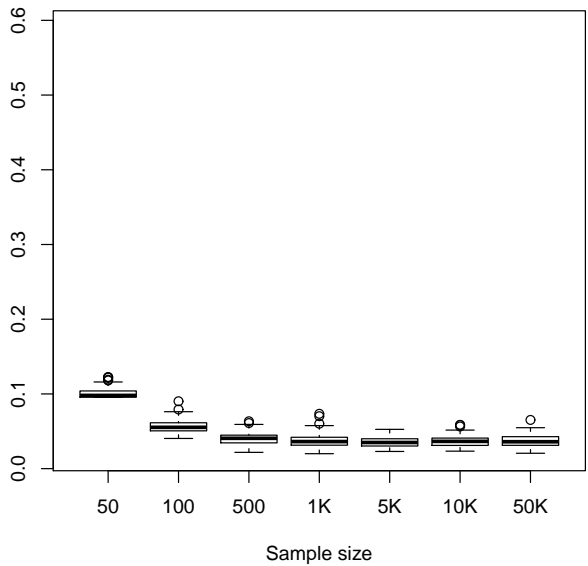}  
\end{center}
\caption{Boxplots of $D_{n, M, M'}$ with $M = 1000$ and $M' =5000$. The sample size $n$ is as indicated and the true convex pmf is $p_0$. See text for details.}
\label{BoxplotAsympConvPi1}
\end{figure}

\begin{figure}[h!]
\begin{center}
\includegraphics[clip=true, width=9cm]{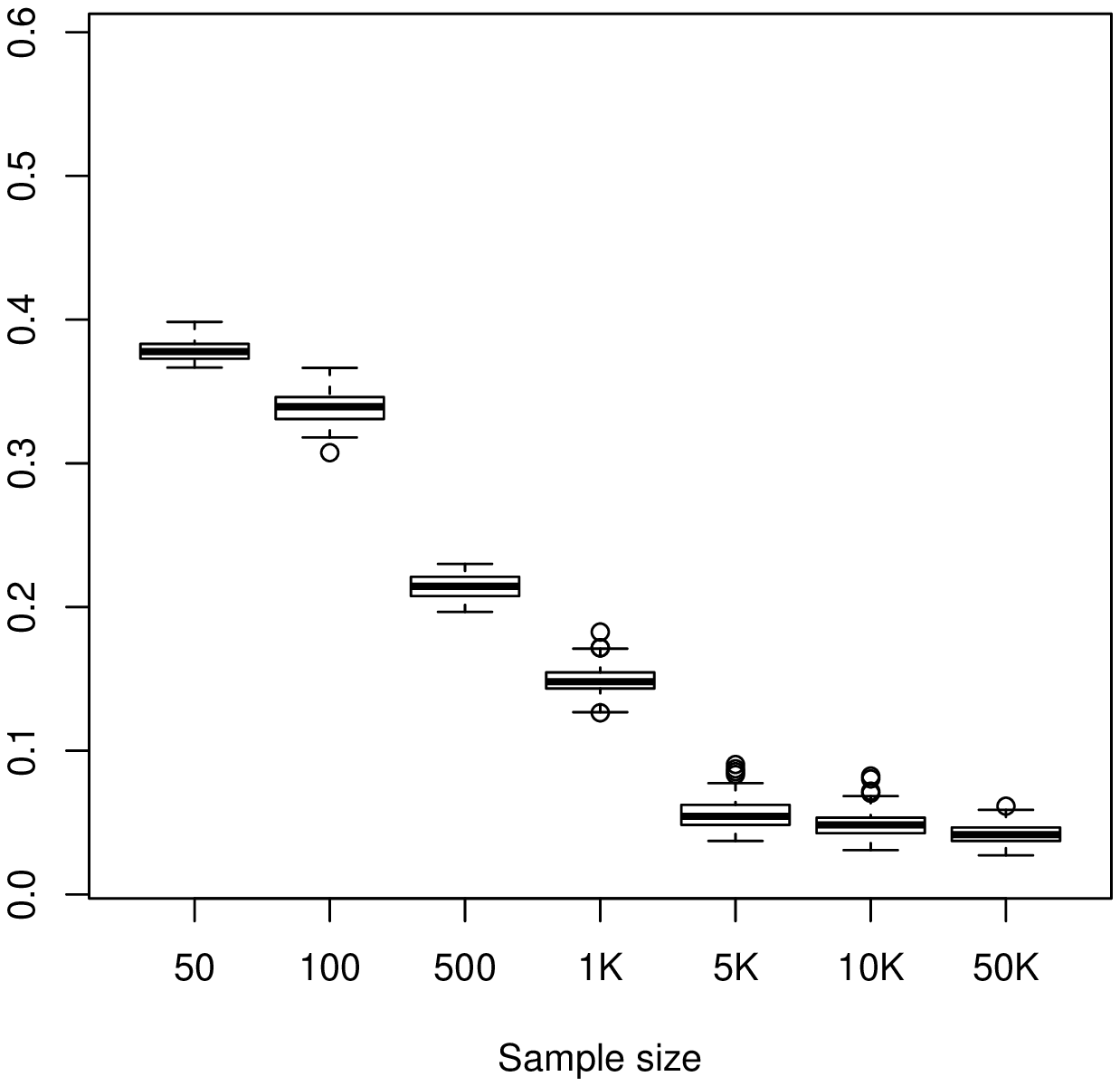}  
\end{center}
\caption{Boxplots of $D_{n, M, M'}$ with $M = 1000$ and $M' =5000$. The sample size $n$ is as indicated and the true convex pmf is $p_1$. See text for details.}
\label{BoxplotAsympConvPi2}
\end{figure}

\begin{figure}[h!]
\begin{center}
\includegraphics[clip=true, width=9cm]{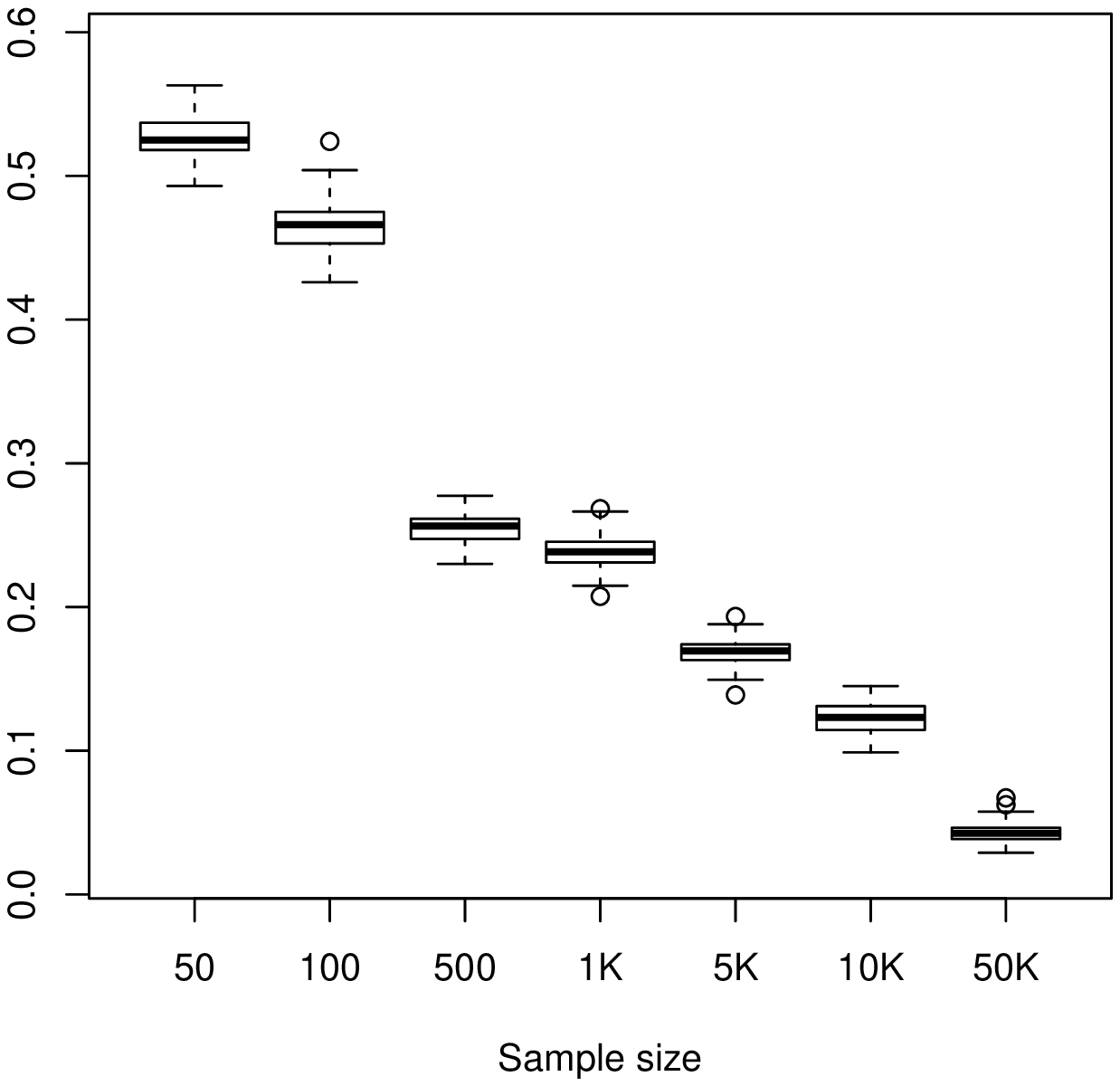}  
\end{center}
\caption{Boxplots of $D_{n, M, M'}$ with $M = 1000$ and $M' =5000$. The sample size $n$ is as indicated and the true convex pmf is $p_2$. See text for details.}
\label{BoxplotAsympConvPi3}
\end{figure}

\begin{figure}[h!]
\begin{center}
\includegraphics[clip=true, width=9cm]{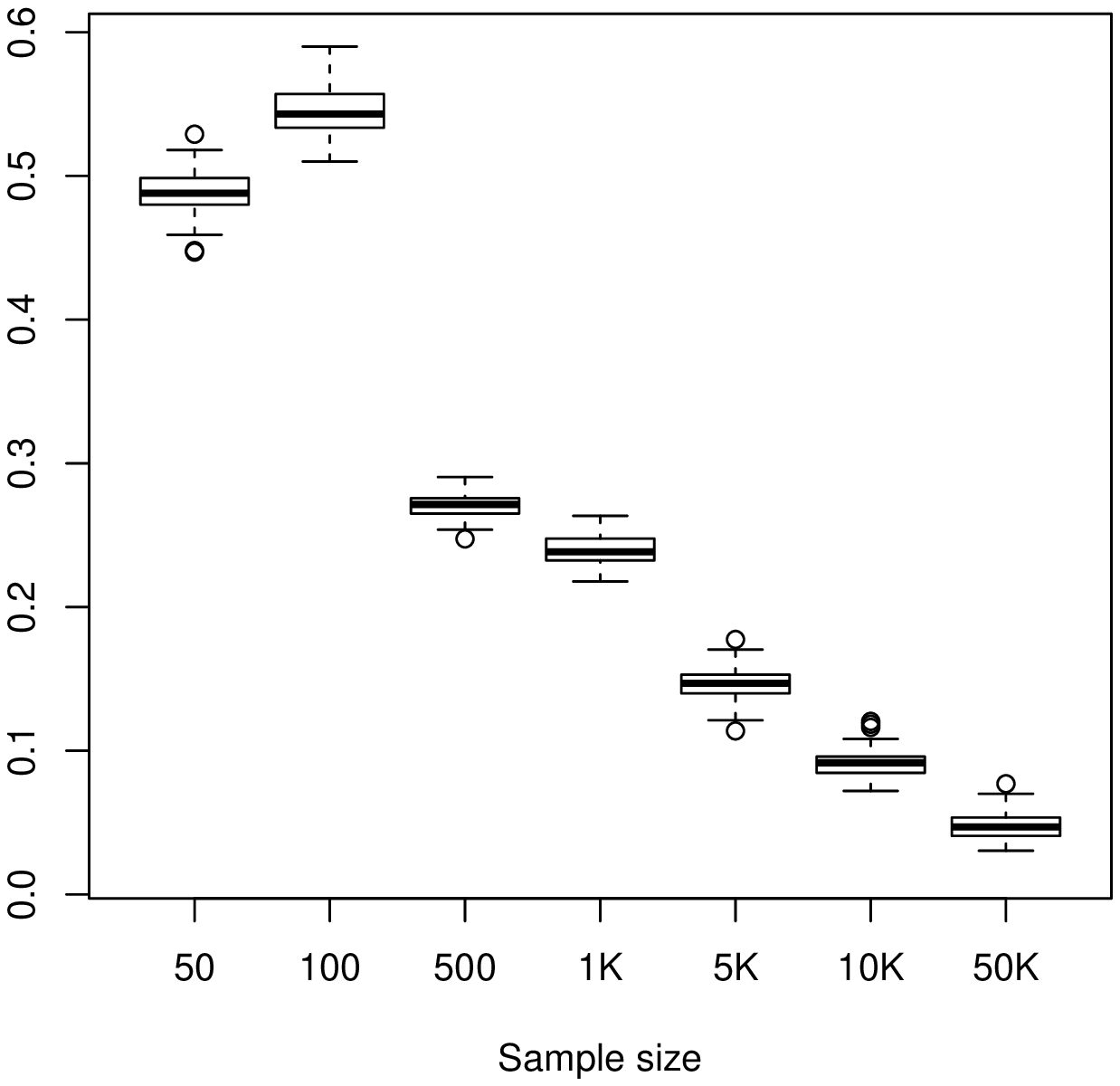}  
\end{center}
\caption{Boxplots of $D_{n, M, M'}$ with $M = 1000$ and $M' =5000$. The sample size $n$ is as indicated and the true convex pmf is $p_3$. See text for details.}
\label{BoxplotAsympConvPi4}
\end{figure}

\begin{figure}[h!]
\begin{center}
\includegraphics[clip=true, width=9cm]{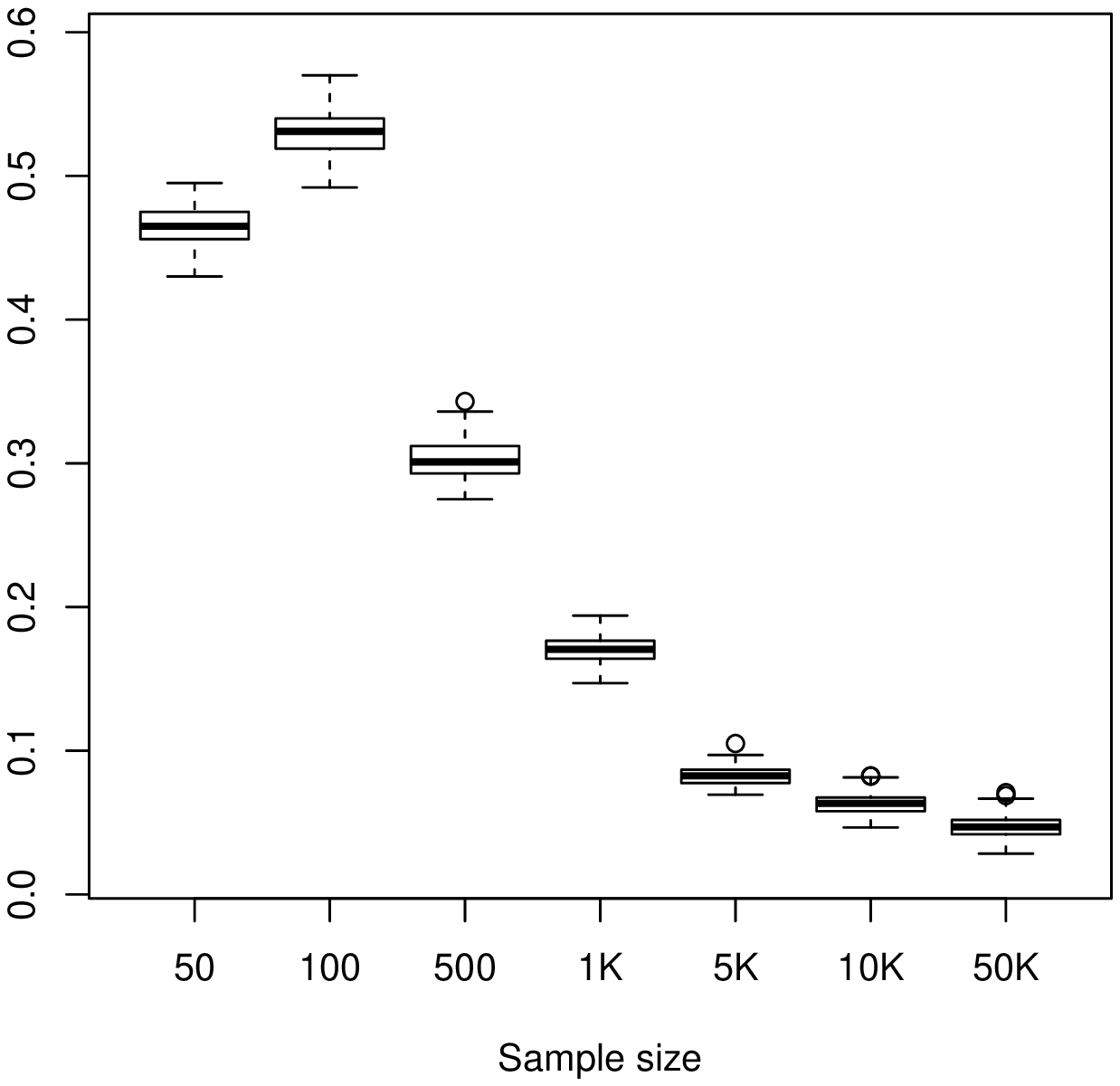}  
\end{center}
\caption{Boxplots of $D_{n, M, M'}$ with $M = 1000$ and $M' =5000$. The sample size $n$ is as indicated and the true convex pmf is $p_4$. See text for details.}
\label{BoxplotAsympConvPi5}
\end{figure}

\begin{figure}[h!]
\begin{center}
\includegraphics[clip=true, width=9cm]{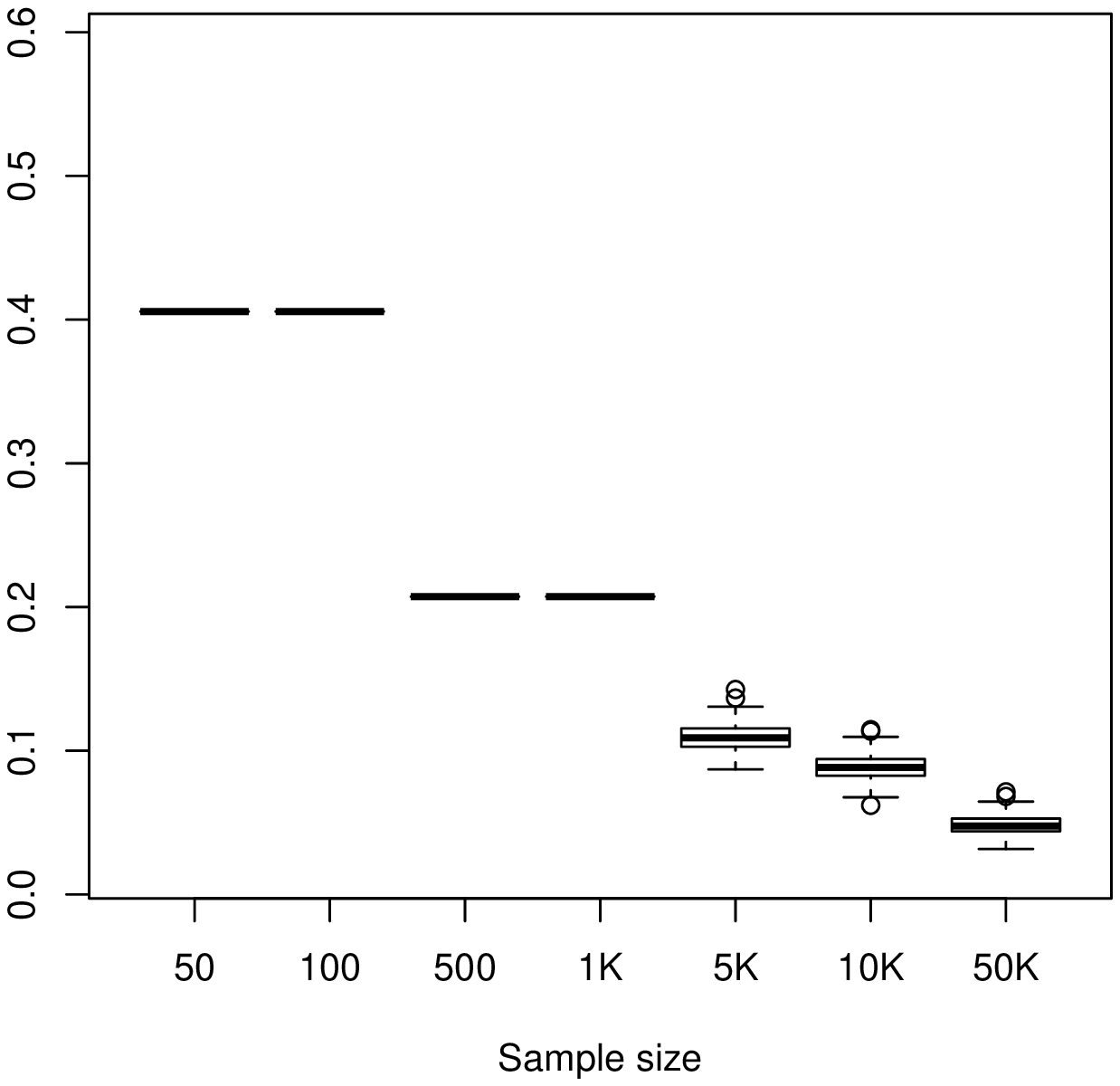}  
\end{center}
\caption{Boxplots of $D_{n, M, M'}$ with $M = 1000$ and $M' =5000$. The sample size $n$ is as indicated and the true convex pmf is the truncated geometric, $p_5$, with success probability equal to $1/2$. See text for details.}
\label{BoxplotAsympConvGeom05}
\end{figure}

To compute the uniform distance between $\widehat{\mathbb{F}}^{(j)}_{n,M}$ and $ \mathbb F^{(j)}_{M'}$, we computed the maximal value of the absolute difference on a discretized grid $\{-3, -2.99, \ldots, 2.99, 3 \}$ with a regular step equal to $0.01$.   The boxplots shown in Figures \ref{BoxplotAsympConvPi1} - \ref{BoxplotAsympConvGeom05} give support to the asymptotic theory of the estimation error of the MLE in case the true pmf is one of the selected convex pmf's $p_0, p_1, \ldots, p_5$. Interestingly, weak convergence seems not to happen at the same speed. For the triangular pmf $p_0$, the boxplots appear to stabilize for $n \ge 1000$ whereas the obtained boxplots for the other distributions seem to indicate that convergence has not been yet attained.  According to our numerical findings in Section \ref{CaptureKnots}, large sample sizes could be required for the estimator to be able to capture these knots. Thus, the slow convergence to the true limit for $p_i, \ 1\ \le i \le 5$ could be partially explained by the fact that the pmf's $p_i, \ 2 \le i \le 5 $ have all interior knots, as opposed the relatively fast convergence in the case of the triangular pmf $p_0$ which has none. We would like to note that all points in $\{1, \ldots, 10 \}$ are interior knots of the truncated geometric pmf $p_5$ since it is strictly convex of its support. Hence, it is the pmf with the largest number of interior knots and also the one for which the convergence seems to be the slowest.  Also, the limit $\widehat g$ reduces  in this particular case to $\WW$. Indeed, the empirical pmf $\hat{p}_{n}$ is convex for large sample sizes (see  Proposition \ref{Knots}) so that $\hat p_{n}=p_{n}$, and hence $\sqrt n(\hat p_{n}-p_{5})$ has the same limit distribution $\WW$ as $\sqrt n(p_{n}-p_{5}).$  Another way of viewing this is to note that the required convexity of the minimizer $\widehat g$ between the knots becomes a superfluous constraint since it is always satisfied by the straight line connecting $\WW$ at two given knots. This is true only in this case because all knots are consecutive. 

We would like to finish this section by adding that it is of course impossible to have a precise statement about the speed of convergence in case the true convex pmf is known to be finitely supported. However, our numerical findings indicate that in applications such as construction of asymptotic confidence bands, one has to keep in mind that moderate sample sizes may not be enough to obtain good coverage. Finally, note that in our assessment we have assumed that the distribution of $\widehat g$ is continuous. We believe this is true but we do not intend to prove it here as it is beyond the scope of this work.

\section{Proofs}\label{sec:proofs}

Before showing Proposition \ref{Consis}, we note that the following inequalities hold true for any sequence $q$. For $r \in [2, \infty]$, we have that 
\begin{eqnarray}\label{Ineqr}
\Vert q \Vert_r  \le  \Vert q \Vert_2, 
\end{eqnarray}
and 
\begin{eqnarray}\label{Ineq2}
\Vert q \Vert_2 \le \Vert q \Vert^{1/2}_\infty \Vert q \Vert^{1/2}_1.
\end{eqnarray}

\medskip
\medskip

\par \noindent \textbf{Proof of Proposition \ref{Consis}.} \  The claim is an immediate consequence of Theorem 4 of \cite{durotetal_13}. Indeed, by the inequality in (\ref{Ineq2}) we have that
\begin{eqnarray*}
\Vert p_n - p_0 \Vert_2 & \le &   \Vert p_n - p_0 \Vert_\infty^{1/2} \Vert  p_n - p_0 \Vert^{1/2}_1 \\
& \le &  \Vert p_n - p_0 \Vert^{1/2}_\infty  \to 0 
\end{eqnarray*}
almost surely by the Glivenko-Cantelli Theorem. Now, by the inequality in (\ref{Ineqr}) and Theorem 4 of \cite{durotetal_13} we can write \begin{eqnarray*}
\Vert \widehat p_n - p_0 \Vert_r & \le &  \Vert \widehat p_n - p_0 \Vert_2 \\
& \le & \Vert p_n - p_0 \Vert_2  \ \to 0, \textrm{a.s.}
\end{eqnarray*} 
 which completes the proof. \hfill $\Box$

\medskip
 \medskip
 
\par \noindent \textbf{Proof of Theorem \ref{BoundProbLSEInt}.}  Theorem 6 of \cite{durotetal_13} implies (\ref{BoundProbLSE}) when $r=\infty$. 
Let $\{0,\dots,S\}$ denote again the support of $p_0$. By Theorem 1 of \cite{durotetal_13}, $\widehat p_{n}$ is a proper pmf. Combining this with Proposition \ref{SuppLSE},  it follows that $\widehat p_n(z) -p_0(z) = 0 $ for $z \ge S+1$. Hence,  for all $z \in \NN$ we have that
\begin{eqnarray*}
\sqrt n \vert F_{\hat p_n}(z) - F_{p_0}(z) \vert  &\le &  \sqrt n \sum_{x =0}^{z \wedge S} \vert \widehat p_n(x) - p_0(x) \vert \\
                                                  & \le  &  (S+1) \sqrt n \Vert \widehat p_n - p_0 \Vert_\infty  = O_p(1).
\end{eqnarray*}
Also, the definition of $H_{\hat p_n}$ and $H_{p_0}$ and  Proposition \ref{SuppLSE} imply that  for all $z \in \NN$
\begin{eqnarray*}
\sqrt n \vert H_{\hat p_n}(z) - H_{p_0}(z) \vert & \le &  \sqrt n \sum_{x=0}^{z \wedge (S+1)-1} \vert F_{\hat p_n}(x) - F_{p_0}(x) \vert \\ 
& \le & (S+1)^2 \sqrt n  \Vert \widehat p_n - p_0 \Vert_\infty = O_p(1) 
\end{eqnarray*}
and the result follows. \hfill $\Box$

\medskip
\medskip

\par \noindent \textbf{Proof of Proposition \ref{SuppLSE}.} \ First, note that the maximal point $\bar s_{n}$ of the support of the empirical pmf $p_{n}$ is $X_{(n)}=\max_{1\leq i\leq n}X_{i}$, and that with probability one, $X_{(n)}=S$ provided that $n$ is sufficiently large. Now, by Theorem 1 of \cite{durotetal_13} we know that $\widehat p_n$ admits a finite support whose maximal point $\widehat s_n \ge \bar s_n$. Therefore, with probability one we have $\widehat s_{n}\geq S$  for $n$ large enough. We show now by contradiction that with probability one there exists $n^*$ such that if $n \ge n^*$ then $\widehat s_n \in \{S, S+1 \}$. Suppose that $\widehat s_n \ge S+2$. By Proposition 1 of \cite{durotetal_13} we know that with probability one there exists $n_0$ such that for $n \ge n_0$, $\widehat p_n$ has to be linear on the set $\{S-1, S, S+1, \ldots, \widehat s_n \}$. But $S+1$ is a knot of $p_0$ which implies by Proposition \ref{Knots} above that with probability one there exits $n^* \ge n_0$ such that for $n \ge n^*$, $S+1$ is also a knot of $\widehat p_n$. This yields a contradiction.  \hfill $\Box$ 

\medskip
\medskip

 \par \noindent \textbf{Proof of Theorem \ref{Exis&Char}.} \ First note that $\cal C(\cal K)$ is a closed convex cone of  $\RR^{S+2}$, so there exists a unique minimizer of $\Phi$ over $\cal C(\cal K)$.

\medskip

Now suppose that $\widehat g$ is the minimizer of $\Phi$ over $\mathcal C(\mathcal K)$. Let $x \in \{1, \ldots, S+2 \}$. Then, for $\epsilon > 0$ the function $k \mapsto \widehat g(k) + \epsilon (x-k)_+ = \widehat{g}(k) + \epsilon T_x(k)$ is clearly in $\mathcal C(\mathcal K)$. Hence,
\begin{eqnarray*}
0 &\le & \lim_{\epsilon \searrow 0} \frac{\Phi(\widehat g  + \epsilon T_x) - \Phi(\widehat g)}{\epsilon} \\
&=& \sum_{k=0}^{S+1} \widehat g(k)   T_x(k) - \sum_{k=0}^{S+1}\WW(k) T_x(k) \\
& = &  \sum_{k=0}^{S+1} \widehat g(k)   (x-k)_+ - \sum_{k=0}^{S+1}\left \{\mathbb U(F_{p_0}(k))-\mathbb U(F_{p_0}(k-1))\right\} (x-k)_+.
\end{eqnarray*}
Recall that $\widehat{\mathbb G}$ is the function defined as
\begin{eqnarray*}
\widehat{\mathbb G}(k) = \sum_{j=0}^k \widehat g(j), \ \textrm{for $k \in \{0, \ldots, S+1 \}$}.
\end{eqnarray*}
For notational convenience, we set $\widehat{\mathbb G}(-1)=0$.
Then, we can rewrite the last inequality as
\begin{eqnarray*}
0 &\le &  \sum_{k=0}^{S+1} \left \{ \widehat{\mathbb G}(k) - \widehat{ \mathbb G}(k-1) \right \} (x-k)_+  - \sum_{k=0}^{S+1} \left \{\mathbb U(F_{p_0}(k))-\mathbb U(F_{p_0}(k-1))\right\} (x-k)_+ \\
  & = & \sum_{k=0}^{S} \widehat{\mathbb G}(k) \left((x-k)_+ - (x-(k+1))_+  \right) + \widehat{\mathbb G}(S+1) \ (x-S-1)_+ \\
  && - \sum_{k=0}^{S} \mathbb U(F_{p_0}(k)) \ \left((x-k)_+ - (x-(k+1))_+  \right)  - \ \mathbb \UU(F_{p_0}(S+1)) \ (x-S-1)_+\\
  &=& \sum_{k=0}^{x-1} \widehat{\mathbb G}(k)  - \sum_{k=0}^{x-1} \mathbb U(F_{p_0}(k)),
\end{eqnarray*}
where the last inequality follows from the fact that 
\begin{eqnarray*}
(x-k)_+ - (x-k -1)_+  = 
\begin{cases}
 1, \ \ \text{if $k \le x-1$} \\
 0, \ \ \text{if $k \ge x$}
\end{cases}
\end{eqnarray*}
together with the fact  that $(x-S-1)_+ =0$ for $x \in \{1, \ldots, S+1\}$, and $(x-S-1)_+ =1$ for $x=S+2$. Thus,
$$
\widehat \HH(x) =  \sum_{k=0}^{x-1} \widehat{\mathbb G}(k)   \ge  \sum_{k=0}^{x-1} \mathbb U(F_{p_0}(k)) = \HH(x)
$$
for all $x \in \{1, \ldots, S+2\}$. Note that at $x=0$, equality of $\widehat \HH$ and $\HH$ is guaranteed by the chosen convention.  Proof of equality in case $x$ is a knot of $\widehat g$  in $\{s_j+1, \ldots, s_{j+1}-1 \}$ uses the fact that the perturbation function $T_x$ satisfies that $\widehat g + \epsilon T_x$ is in $\mathcal C(\mathcal K)$ for $\vert \epsilon \vert $ small enough yielding
$$\lim_{\epsilon \to 0} \frac{1}{\epsilon} (\Phi(\widehat g  +\epsilon T_x) - \Phi(\widehat g )) =0.$$  
We also have equality of $\widehat \HH$ and $\HH$ at the points in $\{s_1, \ldots, s_{m+1}, S+2\}$ since at these points there is no constraint. 

Therefore, if $\widehat g $ is the minimizer of $\Phi$, we have shown that the process $\widehat \HH$ defined on $\{0, \ldots, S+2 \}$ as in (\ref{hatH}) satisfies the inequality and equality in (\ref{CharLSIneq2}). Conversely, suppose that $\widehat g  \in \mathcal C(\mathcal K)$ such that the process $\widehat \HH$ in (\ref{hatH}) satisfies (\ref{CharLSIneq2}). Let $g \in \mathcal C(\mathcal K)$. We will show now that $\Phi(g) \ge \Phi(\widehat g)$. We have
\begin{eqnarray*}
\Phi(g) - \Phi(\widehat g ) 
&=&
\frac{1}{2}\sum_{k=0}^{S+1}(g(k)-\widehat g(k))^2+\sum_{k=0}^{S+1}(g(k)-\widehat g(k))(\widehat g(k)-\WW(k))\\
& \ge & \sum_{k=0}^{S+1}  (g(k) - \widehat g (k)) \ \widehat g(k)  \\
&& \ - \sum_{k=0}^{S+1} \left(g(k) - \widehat g (k) \right) \left (\mathbb U(F_{p_0}(k))-\mathbb U(F_{p_0}(k-1))\right) \\
& = & \sum_{k=0}^{S+1} (g(k) - \widehat g  (k)) (\widehat{\mathbb D}(k) - \widehat{\mathbb D}(k-1)), 
\end{eqnarray*}
where $\widehat{\mathbb D}(k) = \widehat{\mathbb G}(k) - \mathbb U(F_{p_0}(k))$.  Now, similar as above, for all $x\in\{1,\dots,S+2\}$ we have
\begin{eqnarray}\label{sumD}
\sum_{k=0}^{S+1} (\widehat{\mathbb D}(k) - \widehat{\mathbb D}(k-1)) (x-k)_+ &=& 
 \sum_{k=0}^{x-1} \widehat{\mathbb D}(k) \notag\\
 & =& \widehat{\mathbb H}(x) - \HH(x) \notag\\
 &\geq& 0,
\end{eqnarray}
with equality if $x$ is a knot of $\widehat g $ or $x \in \{s_0, \ldots, s_{m+1}, S+2\}$. To conclude,
we will use the fact that an arbitrary element $g \in \mathcal C(\mathcal K)$ can be written as
 \begin{equation}\label{RepForm2}
g(k) =  \alpha +\sum_{j=1}^{m+1}c_j (s_j-k)_+ + \sum_{j=1}^{m+1}\sum_{i=1}^{J_j} c_{j,i} (z_{j,i} - k)_+\\
\end{equation}
 for all $k=0,\dots,S+1$, with $z_{j,1},\dots ,z_{j,J_{j}}$ the knots of $g$ in $\{s_{j-1}+1,\dots,s_{j}-1\}$ for $j=1,\dots,m+1$, and where $\alpha, c_1, c_2, \ldots, c_{m+1} $ are real numbers, and $c_{j, i} > 0$ for $j=1,\dots,m+1$ and $i=1,\dots,J_{j}$. This comes from the fact that  any finite convex sequence $p=\{p(0),\dots,p(K)\}$ for some $K>0$, admits the (spline) representation
\begin{eqnarray}\label{RepForm}
p(k) = a + \gamma_1 (s_1 - k)_+ + \ldots + \gamma_p (s_p - k)_+ +\gamma_{p+1}(K-k)_{+}
\end{eqnarray}
where $a$ is a real number, $\gamma_i > 0$ and $0< s_1 < \ldots < s_p < K $ are the interior knots of $p$.  Using the spline representation in (\ref{RepForm2})  together with \eqref{sumD}, where we recall that we have an equality for $x=s_{1},\dots,s_{m+1}$, it follows that
\begin{eqnarray*}
&&\sum_{k=0}^{S+1} (\widehat{\mathbb D}(k) - \widehat{\mathbb D}(k-1)) g(k) \\
&&\qquad\qquad =\alpha\widehat{\mathbb D}(S+1) +\sum_{j=1}^{m+1}\sum_{i=1}^{J_j} c_{j,i} \sum_{k=0}^{S+1} (\widehat{\mathbb D}(k) - \widehat{\mathbb D}(k-1)) (z_{j,i}-k)_+ \\
&&\qquad\qquad \geq\alpha\widehat{\mathbb D}(S+1),
\end{eqnarray*}
where in the last inequality we used the fact that $c_{j,i}\geq 0$ for all $j,i$. Now,  the boundary conditions $\widehat \HH(S+1) = \HH(S+1)$ and $\widehat \HH(S+2) = \HH(S+2)$ in (\ref{CharLSIneq2}) imply that 
\begin{eqnarray*}
\widehat{\mathbb D}(S+1)&=& \sum_{k=0}^{S+1} \widehat g(k)\\
  &=&  \widehat \HH(S+2) - \widehat \HH(S+1) \\ 
& = & \HH(S+2) - \HH(S+1) = \UU(F_{p_0}(S+1)) = 0,
\end{eqnarray*}
since $F_{0}(S+1)=1$ and $\UU$ is a standard Brownian bridge. We arrive at
 $$
\sum_{k=0}^{S+1} (\widehat{\mathbb D}(k) - \widehat{\mathbb D}(k-1))  g (k) \geq 0
$$
and similarly, since we have an equality in \eqref{sumD} if $x$ is a knot of $\widehat g$,
$$
\sum_{k=0}^{S+1} (\widehat{\mathbb D}(k) - \widehat{\mathbb D}(k-1)) \widehat g (k) = 0.
$$
It follows that $\Phi(g) \ge  \Phi(\widehat g )$ and that $\widehat g $ is the minimizer of $\Phi$.  \hfill $\Box$

\medskip
\medskip

\par \noindent \textbf{Proof of Theorem \ref{Asymp}.}  \ For $z \in \NN$ define  
\begin{eqnarray*}
\YY_n(z) =: \sum_{k=0}^{z-1} \sqrt n \left(F_{p_n}(k) - F_{p_0}(k)\right) = \sum_{k=0}^{z-1} \UU_n(F_{p_0}(k)), 
\end{eqnarray*}
\begin{eqnarray*}
\widehat Y_n(z) =:  \sum_{k=0}^{z-1} \sqrt n \left(F_{\hat p_n}(k) - F_{p_0}(k)\right),
\end{eqnarray*}
and recall that
\begin{eqnarray*}
\HH(z) = \sum_{x=0}^{z-1} \UU(F_{p_0}(k))
\end{eqnarray*}
with $\YY_n(0) = \widehat Y_n(0) = \HH(0) =0$. 

\medskip

It follows from the characterization of $\widehat p_n$ in (\ref{Fenchel}) that
\begin{eqnarray*}
\widehat Y_n(x)   
\begin{cases}
\ge  \YY_n(x) ,  \  \text{for $x \ \in \{0, \ldots, S+2 \}$}  \\
= \YY_n(x), \  \text{if $x  \ \in \{0, \ldots, S+2 \}$ is a knot of $\widehat p_n$}.
\end{cases}
\end{eqnarray*}
By standard results on weak convergence of empirical processes we have that
\begin{eqnarray*}
\YY_n(x) \to_d  \HH(x), \ \text{for $x \ \in \{0, \ldots, S+2 \}$}.  
\end{eqnarray*}
Now, Theorem \ref{BoundProbLSEInt} implies that there exists a subsequence $\{\widehat Y_{n'}\}_{n'}$ which weakly converges to some $\LL$  on $\{0, \ldots, S+2\}$. In what follows, we will use the Skorokhod representation to assume that convergences of $\{\widehat Y_{n'}\}_{n'}$ and $\{\YY_{n'}\}_{n'}$ to their respective limits happen almost surely. The goal now is to show that $\LL$ and $\widehat\HH$, defined in Theorem \ref{Exis&Char} above, are equal with probability one. Let us define
$$
\widetilde{g}(x) := \LL(x+1) + \LL(x-1) - 2 \LL(x)
$$
for $x \in \{0, \ldots, S+1 \}$ with  $\LL(-1)=0$. Note that we have that 
\begin{eqnarray}\label{IdS+2}
\LL(S+2)= \HH(S+2)= \sum_{k=0}^{S+1} \UU(F_{p_0}(k)) = \sum_{k=0}^{S} \UU(F_{p_0}(k)).
\end{eqnarray}
Indeed, by Fubini's theorem we have that 
\begin{eqnarray*}\label{Th5}
\sum_{k =0}^{S+1} (1-F_{\hat p_n}(k)) & = & \sum_{k=0}^{S+1} k \widehat p_n(k) \notag \\
                                                        & =  & \sum_{k=0}^{S+1} k p_n(k), \ \textrm{by Theorem 5 of \cite{durotetal_13}}  \notag\\
                                                        &= & \sum_{k=0}^{S+1} (1- F_{p_n}(k)). \notag
\end{eqnarray*}
This in turn implies that
\begin{eqnarray*}
\widehat Y_n(S+2) & = & \sum_{k=0}^{S+1} \sqrt n (F_{\hat p_n}(k) - F_{p_0}(k)) \\
                         & = &  \sum_{k=0}^{S+1} \sqrt n (F_{p_n}(k) - F_{p_0}(k)) =  \YY_n(S+2) (= \YY_n(S+1)) 
\end{eqnarray*}
and the claim follows by passing $n' \to \infty$.  Hence,  with probability one we have that
\begin{eqnarray*}
\LL(x)   
\begin{cases}
\ge  \HH(x) , &  \  \textrm{for $x \ \in \{0, \ldots, S+2 \}$}  \\
= \HH(x), & \  \textrm{if  $x \in \{s_0, \ldots, s_{m+1}, S+2 \} $ } \\
  & \hspace{0.2cm}  \textrm{or $x$ is a knot of $\widetilde{g}$ in $\{s_j+1, \ldots, s_{j+1}- 1 \}$} \\
  & \hspace{0.2cm} \textrm{for some $j \in \{0, \ldots, m \}$}.
\end{cases}
\end{eqnarray*}
Equality of $\LL$ and $\HH$ at points in $\{0,s_1, \ldots, s_{m+1}, S+2 \}$ in the last assertion follows from the chosen convention at $s_0 =0$, Proposition \ref{Knots},  the equalities $\widehat Y_{n'}(s_j) =  \YY_{n'}(s_j) $ for $j=0, \ldots, m+1$ with probability one and $n'$ large enough and the identity in (\ref{IdS+2}). Equality of $\LL$ and $\HH$ at the knots of $\widetilde g$ in $\{s_j+1, \ldots, s_{j+1}- 1 \}$ follows from the fact that $x$ is a knot of $\widehat p_{n'}$ in $\{s_j+1, \ldots, s_{j+1}- 1 \}$ if and only if it is a knot of $\sqrt {n'} (\widehat p_{n'} - p_0)$ using linearity of $p_0$ between two successive knots. Thus, if $x$ is a knot of $\widetilde{g}$, then as soon as $n'$ is large enough $x$ is also a knot of $$z \mapsto \widehat Y_{n'}(z+1) + \widehat Y_{n'}(z-1) - 2 \widehat Y_{n'}(z) = \sqrt {n'} (\widehat p_{n'}(z) - p_0(z)).$$ This in turn implies that $\widehat Y_{n'}(x) = \YY_{n'}(x)$ implying after passing to the limit that $\LL(x) = \HH(x)$ almost surely.   Furthermore,  $\widetilde{g}$ is clearly in $\mathcal C(\mathcal K)$.

Therefore, it follows  from Theorem \ref{Exis&Char} that $\widetilde{g}$ must be equal to the minimizer of $\Phi$ defined in (\ref{Phi}). Thus,  there exists a version of $\widehat g$ such that
\begin{eqnarray*}
\widehat g(x) = \LL(x+1) + \LL(x-1) - 2 \LL(x)  = \widehat{\HH}(x+1) + \widehat{\HH}(x-1) - 2 \widehat{\HH}(x) 
\end{eqnarray*}
for $x \in \{0, \ldots, S+1 \}$. Put $\Delta = \widehat{\HH} - \LL$. Then, for $x \in \{0, \ldots, S+1 \}$, we have that  $\Delta(x+1) - \Delta(x) = \Delta(x) - \Delta(x-1)$.  Thus, $\Delta(x) - \Delta(x-1) = \Delta(0) - \Delta(-1) = 0$ since $\LL(-1) = \widehat{\HH}(-1) = 0$ and $\LL(0) = \widehat{\HH}(0) =0$.  We conclude that $\Delta(x) = \Delta(S+1) = 0$ since $\widehat{\HH}(S+1) = \LL(S+1) = \HH(S+1)$.  Since we also have $\widehat{\HH}(S+2) = \LL(S+2) = \HH(S+2)$, it follows that $\widehat{\HH} = \LL$ on $\{0, \ldots, S+2 \}$. Now, from an arbitrary subsequence $\widehat{Y}_{n'}$ we can extract a further subsequence $\{\widehat Y_{n''}\}_{n''}$ converging to $\widehat{\HH}$. Since the limit is the same for any such subsequence, we conclude that $\{\widehat Y_n\}_n$ converges weakly to $\widehat{\HH}$ on $\{0, \ldots, S+2\}$ with probability one.  This in turn implies that the following converges
\begin{eqnarray*}
\sqrt n (H_{\hat p_n} - H_{p_0})  \Rightarrow \widehat{\HH}
\end{eqnarray*} 
\begin{eqnarray*}
\sqrt n (F_{\hat p_n} - F_{p_0})  \Rightarrow \widehat{\mathbb G}
\end{eqnarray*} 
and
\begin{eqnarray*}
\sqrt n (\widehat p_n - p_0)  \Rightarrow  \widehat g
\end{eqnarray*} 
occur jointly, and the proof is complete. \hfill $\Box$

\medskip
\medskip

\par \noindent \textbf{Proof of Theorem \ref{LocalLeft}.}  \ Existence and uniqueness of the minimizer both follow from the projection theorem on closed convex cones in $\RR^{s+1}$. With similar arguments as for the proof of Theorem \ref{Exis&Char}, it can be shown that an arbitrary element $\widehat g^{\leq s}\in{\cal C}^{\leq s}(\cal K)$ is the minimizer if and only if the process $ \widehat{\HH}^{\leq s}$ defined on $\{0,\ldots, s+1\}$ by 
\begin{eqnarray*}
\widehat\HH^{\leq s}(x) = 
\begin{cases}
\sum_{k=0}^{x-1} \sum_{j=0}^k \widehat g^{\leq s}(j), \ \text{if $x \in \{1, \ldots, s+1\}$} \\
 0, \ \hspace{1.25cm} \text{if $x=0$}
\end{cases}
\end{eqnarray*}
satisfies
$$
\widehat \HH^{\leq s}(x)
\ge  \HH(x)$$ 
for all $x \in \{0, \ldots, s+1\}$, 
with an equality 
if  $x \in \{s_{0},\dots,s_{m+1},s+1\}$ or $x$ is a knot of $\widehat g^{\leq s}$ in $\{s_{j}+1,\dots,s_{j+1}-1\}$ for some $j=0,\dots,m$ with $s_{j+1}<s$.
Consider the restriction $\widehat g^{\leq s}=(\widehat g(0),\dots,\widehat g(s))$ of $\widehat g$ to $\{0,\dots,s\}$. A point $x\in\{1,\dots,s-1\}$ is a knot of $\widehat g^{\leq s}$ if, and only if, it is a knot of $\widehat g$. Therefore, it immediately follows from the characterization of $\widehat g$ given in Theorem \ref{Exis&Char} that the minimizer of $\Phi^{\leq s}$ over ${\mathcal C}^{\leq s}(\mathcal K)$ is equal to $(\widehat g(0),\dots,\widehat g(s))$ if, and only if, $\widehat\HH(s+1)=\HH(s+1)$. Since we already have $\widehat\HH(s)=\HH(s)$, this is equivalent to 
$$\widehat\HH(s+1)-\widehat\HH(s)=\HH(s+1)-\HH(s),$$
that is \eqref{NSC Left}.

\medskip

To prove the last assertion, note that from Theorem \ref{Asymp}, it follows that $\sqrt n(\widehat p_{n}(0)-p_{0}(0),\dots,\widehat p_{n}(s)-p_{0}(s))$ converges in distribution to $(\widehat g(0),\dots,\widehat g(s))$ as $n\to\infty$. Thus, it converges in distribution to the minimizer of $\Phi^{\leq s}$ over ${\mathcal C}^{\leq s}(\mathcal K)$ if, and only if, 
\eqref{NSC Left} holds true with probability one. On the other hand, fromTheorem \ref{Asymp}, one also has
\begin{equation}\label{cvJointe}
\sqrt n\begin{pmatrix}F_{\hat p_n}-F_{p_0}\\ F_{p_n}-F_{p_0}\end{pmatrix} \Rightarrow \begin{pmatrix}
\widehat{\mathbb G}\\ \UU_{n}(F_{p_{0}})\end{pmatrix}\end{equation}
as $n\to\infty$. Therefore, one can have \eqref{NSC Left} with probability one if, and only if, 
$$\sqrt n\Big((F_{\hat p_n}(s)-F_{p_{0}}(s))-(F_{p_n}(s)-F_{p_{0}}(s))\Big)$$
converges in probability to zero as $n\to\infty$. This is equivalent to  \eqref{NSC Left n}, so the proof of the theorem is complete.  \hfill $\Box$

\medskip
\medskip

\par \noindent \textbf{Proof of Theorem \ref{LocalRight}.} \ Existence and uniqueness of the minimizer both follow from the projection theorem on closed convex cones in $\RR^{S-s}$. With similar arguments as for the proof of Theorem \ref{Exis&Char}, it can be shown that an arbitrary element $\widehat g^{\geq s}\in{\cal C}^{\geq s}(\cal K)$ is the minimizer if and only if the process $ \widehat{\HH}^{\geq s}$ defined on $\{s,\ldots, S+2\}$ by 
\begin{eqnarray*}
\widehat\HH^{\geq s}(x) = 
\begin{cases}
\sum_{k=s}^{x-1} \sum_{j=s}^k \widehat g^{\geq s}(j), \ \text{if $x \in \{s+1, \ldots, S+2\}$} \\
0, \ \hspace{1.25cm} \text{if $x=s$}
\end{cases}
\end{eqnarray*}
satisfies
$$
\widehat \HH^{\geq s}(x)
\ge  \HH^{\geq s}(x)$$
for all $x \in \{s, \ldots, S+2\}$, with an equality if  $x \in \{s_{0},\dots,s_{m+1},S+2\}$ or $x$ is a knot of $\widehat g^{\geq s}$ in $\{s_{j}+1,\dots,s_{j+1}-1\}$ for some $j=0,\dots,m$ with $s_{j}\geq s$,
where $\HH^{\geq s}(s)=0$ and
$$\HH^{\geq s}(x)=\sum_{k=s}^{x-1} \sum_{j=s}^k \WW(j)$$
for all $x\in\{s+1,\dots,S+2\}$. The connection between $\HH$ and $\HH^{\geq s}$ is as follows. For all $x\in\{s,\dots,S+2\}$ one has
\begin{eqnarray*}
\HH(x)&=&\sum_{k=0}^{x-1} \sum_{j=0}^k \WW(j)\\
&=&\sum_{k=0}^{s-1} \sum_{j=0}^k \WW(j)+\sum_{k=s}^{x-1} \sum_{j=0}^{s-1} \WW(j)+\HH^{\geq s}(x)\\
&=&\HH(s)+(x-s)\Big(\HH(s)-\HH(s-1)\Big)+\HH^{\geq s}(x).
\end{eqnarray*}
Similarly,
$$
\widehat \HH(x)=\widehat \HH(s)+(x-s)\Big(\widehat \HH(s)-\widehat \HH(s-1)\Big)+\widehat \HH^{\geq s}(x)
$$
for all $x\in\{s,\dots,S+2\}$, so that
\begin{equation}\label{connect}
\widehat \HH^{\geq s}(x)-\HH^{\geq s}(x)=\widehat \HH(x)-\HH(x)+(x-s)\Big(\widehat \HH(s-1)-\HH(s-1)\Big),
\end{equation}
using that $\widehat\HH(s)=\HH(s)$.

Consider  $\widehat g^{\geq s}=(\widehat g(s),\dots,\widehat g(S+1))$. A point $x\in\{s+1,\dots,S\}$ is a knot of $\widehat g^{\geq s}$ if, and only if, it is a knot of $\widehat g$. Therefore, it immediately follows from the characterization of $\widehat g$ given in Theorem \ref{Exis&Char}, together with \eqref{connect},  that the minimizer of $\Phi^{\geq s}$ over ${\mathcal C}^{\geq s}(\mathcal K)$ is equal to $(\widehat g(s),\dots,\widehat g(S+1))$ if, and only if, $\widehat\HH(s-1)=\HH(s-1)$. Since we already have $\widehat\HH(s)=\HH(s)$, this is equivalent to 
$$\widehat\HH(s)-\widehat\HH(s-1)=\HH(s)-\HH(s-1),$$
that is \eqref{NSC Right}.

\medskip

To prove the last assertion, note that from Theorem \ref{Asymp}, it follows that $\sqrt n(\widehat p_{n}(s)-p_{0}(s),\dots,\widehat p_{n}(S+1)-p_{0}(S+1))$ converges in distribution to $(\widehat g(s),\dots,\widehat g(S+1))$ as $n\to\infty$. Thus, it converges in distribution to the minimizer of $\Phi^{\geq s}$ over ${\mathcal C}^{\geq s}(\mathcal K)$ if, and only if, 
\eqref{NSC Right} holds true with probability one. According to \eqref{cvJointe}, this happens if, and only if, 
$$\sqrt n\Big((F_{\hat p_n}(s-1)-F_{p_{0}}(s-1))-(F_{p_n}(s-1)-F_{p_{0}}(s-1))\Big)$$
converges in probability to zero as $n\to\infty$. This is equivalent to  \eqref{NSC Right n}, so the proof of the theorem is complete.  \hfill $\Box$

\medskip

\par \noindent{Characterizing convexity for a truncated geometric pmf.} \  A truncated pmf on $\{0, \ldots, S \}$ with success probability $\alpha \in [0,1]$ is given by 
\begin{eqnarray*}
q_\alpha(i) = \left \{
\begin{array}{ll}
\frac{\alpha}{1 - (1- \alpha)^{S+1}}  (1-\alpha)^i, \ \ i \in \{0, \ldots, S\} \\
0 \ \ \hspace{1cm} \ \textrm{otherwise}.
\end{array}
\right.
\end{eqnarray*} 
Hence, $q_\alpha$ is convex on $\NN$ if and only if 
\begin{eqnarray*}
q_\alpha(S+1) - q_{\alpha}(S) \ge q_\alpha(S) - q_{\alpha}(S-1)
\end{eqnarray*} 
or equivalently
\begin{eqnarray*}
-(1-\alpha)^S  \ge (1-\alpha)^{S} - (1-\alpha)^{S-1} 
\end{eqnarray*}
that is $\alpha \ge 1/2$. \hfill $\Box$

\medskip
\medskip

\par \noindent \textbf{Proof of Lemma \ref{lem:geometric}.} \ We have $\Delta p(k)\geq 0$ for all $k<S$ by convexity of the function $k\mapsto q^k$. We also have $\Delta p(k)\geq 0$ for all integers $k>S+1$ since for those $k$, $\Delta p(k)= 0$. Now, $\Delta p(S+1)=p(S)\geq 0$, and
$$\Delta p(S)=\frac{q^{S-1}(1-q)}{1-q^{S+1}}(1-2q),$$
which is $\geq 0$ if, and only if, $q\leq 1/2$. We conclude that $\Delta p(k)\geq 0$ for all $k\in\mathbb{N}$ if, and only if, $q\leq 1/2$. \hfill $\Box$ 

\paragraph{Acknowledgments} \  The authors are very grateful to Sylvie Huet for many helpful comments, stimulating discussions on the subject, and also for a careful reading of the paper. 

\bibliographystyle{ims}
\bibliography{PaperAsymptLSE_ArXiv}

\end{document}